\documentclass[12pt]{article}

\newlength{\spse}
\usepackage{latexsym,amsmath,amssymb,amsfonts,amscd,multirow}
\usepackage{epsfig,verbatim,epstopdf,graphics}
\usepackage{bm,color}

\newtheorem{thm}{Theorem}[section]

\newtheorem{lem}[thm]{Lemma}
\newtheorem{exa}[thm]{Example}
\newtheorem{rem}[thm]{Remark}
\newtheorem{prop}[thm]{Property}

\newcommand{\sign}{\operatorname{sign}}

\newcommand{\bk}{\mathbf{k}}
\newcommand{\bx}{\mathbf{x}}

\newfont{\iams}{msbm9}

\newcommand{\commentbis}[1]{}
\newcommand{\be}{\begin{eqnarray}}
\newcommand{\ee}{\end{eqnarray}}
\newcommand{\beno}{\begin{eqnarray*}}
\newcommand{\eeno}{\end{eqnarray*}}

\newcommand{\barr}[1]{\begin{array}{#1}}
\newcommand{\earr}{\end{array}}

\newcommand{\beq}{\begin{equation}}
\newcommand{\eeq}{\end{equation}}
\newcommand{\beqa}{\begin{eqnarray}}
\newcommand{\eeqa}{\end{eqnarray}}

\newcommand{\bv}{{\bf v}}

\newcommand{\bV}{{\bf V}}
\newcommand{\bn}{{\bf n}}

\newcommand{\bP}{{\bf P}}

\newcommand{\bzero}{\mathbf{0}}
\newcommand{\bone}{\mathbf{1}}
\newcommand{\bl}{\mathbf{l}}
\newcommand{\bi}{\mathbf{i}}

\newcommand{\bj}{\mathbf{j}}
\newcommand{\bW}{\mathbf{W}}

\newcommand{\bK}{\mathbf{K}}
\newcommand{\ba}{{\bm{\alpha}}}
\newcommand{\bb}{{\bm{\beta}}}
\newcommand{\bq}{\mathbf{q}}

\makeatletter
\newcommand{\Rmnum}[1]{\expandafter\@slowromancap\romannumeral #1@}
\makeatother

\title
{Sparse Grid Discontinuous Galerkin Methods for High-Dimensional
Elliptic Equations}

\author{
Zixuan Wang \thanks{Department of Mathematics, Michigan State University, East Lansing, MI 48824 U.S.A.
{\tt wangzix1@msu.edu}}
\and
Qi Tang \thanks{Department of Mathematical Sciences, Rensselaer Polytechnic Institute, Troy, NY 12180 U.S.A.
{\tt tangq3@rpi.edu}}
\and
 Wei Guo
\thanks{Department of Mathematics, Michigan State University,
East Lansing, MI 48824 U.S.A.
 {\tt wguo@math.msu.edu}}
\and
 Yingda Cheng
\thanks{Department of Mathematics, Michigan State University,
East Lansing, MI 48824 U.S.A.
 {\tt ycheng@math.msu.edu}. Research is supported by NSF grant DMS-1318186.}
}

\date{\today}

\begin{document}

\maketitle

\begin{abstract}
This paper constitutes our initial effort in developing sparse grid discontinuous Galerkin (DG) methods for high-dimensional partial differential equations (PDEs).   Over the past few decades,  DG methods have gained popularity in many applications due to their distinctive features.
However, they are often deemed too costly because of the large number of degrees of freedom of the approximation space, which are the main bottleneck for simulations in high dimensions. In this paper, we develop sparse grid DG methods for   elliptic equations with the aim of breaking the \emph{curse of dimensionality}. Using a hierarchical basis representation, we construct a sparse finite element approximation space, reducing the degrees of freedom  from the standard {$O(h^{-d})$ to $O(h^{-1}|\log_2 h|^{d-1})$} for $d$-dimensional problems, where $h$ is the uniform mesh size in each dimension. Our method, based on the interior penalty (IP) DG framework, can achieve  accuracy of  $O(h^{k}|\log_2 h|^{d-1})$ in the energy norm, where $k$ is the degree of   polynomials used.
Error estimates  are provided and confirmed by numerical tests in multi-dimensions.

\end{abstract}

{\bf Keywords:} discontinuous Galerkin methods;   interior penalty methods; sparse grid; high-dimensional partial differential equations.

\section{Introduction}

Elliptic equations have been used successfully in many areas of modeling, such as aeroacoustics, electro-magnetism, oil recovery simulation, weather forecasting, etc. Their numerical treatment includes the standard finite difference methods, collocation methods, Galerkin finite element methods, least squares methods, among many others. Discontinuous Galerkin (DG) methods work well for purely hyperbolic problems due to the discontinuous nature of the solutions. Yet, they are also proven to be effective for elliptic equations. Since the 1970s, the interior penalty (IP) DG  methods have been proposed \cite{baker1977finite, douglas1976interior, wheeler1978elliptic, arnold1982interior} and  generalized \cite{oden1998discontinuoushpfinite, riviere2001priori, dawson2004compatible}.   A  review of DG methods for elliptic equations including discussions about other DG methods, such as the method of Bassi and Rebay \cite{Bassi_1997_JCP_DFEM} and  the local DG (LDG) method \cite{cockburn1998local, castillo2000priori, cockburn2001superconvergence}, can be found in \cite{Arnold_2002_SIAM_DG}.
It is generally understood that the advantages of DG methods lie in their flexibility in choosing the discretization meshes and the approximation spaces. However, the methods are often deemed too costly due to the  large degrees of freedom of the approximation space. Such drawback is more prominent when dealing with high-dimensional equations arising from real-world applications,  such  as kinetic simulations, stochastic analysis, and mathematical modeling in finance or statistics.

This paper concerns the computation of high-dimensional elliptic equations, for which the main challenge  is commonly known as the \emph{curse of dimensionality} \cite{bellman1961adaptive}. This term refers to the fact that the computational cost and storage requirements   scale as  $O(h^{-d})$ for a $d$-dimensional problem, where   $h$ denotes the mesh size in one coordinate direction. This challenge   typically can not be resolved through barely increasing computational resources, and it requires the improvement of numerical techniques as well as  efficient computational implementations.
The sparse grid techniques \cite{bungartz2004sparse, garcke2013sparse},  introduced by Zenger \cite{zenger1991sparse}, have been developed as a major tool to break the curse of dimensionality of grid-based approaches. The idea relies on a tensor product hierarchical basis representation, which can reduce the degrees of freedom   from   {$O(h^{-d})$ to $O(h^{-1}|\log_{2}h|^{d-1})$ for $d$-dimensional problems without compromising much accuracy.} The fundamentals of sparse grid techniques can be traced back to Smolyak \cite{smolyak1963quadrature} for numerical integration, and they are closely related to 
hyperbolic cross \cite{babenko1960approximation,temlyakov1986approximations}, boolean method \cite{delvos1982d}, discrete blending method \cite{baszenski1992blending}, and splitting extrapolation method \cite{liem1995splitting} in the literature. The construction of the scheme is based on  balancing the cost complexities and the accuracy of the scheme by seeking a proper truncation of the tensor product hierarchical bases, which can be formally derived by solving an optimization problem of cost/benefit ratios  \cite{griebel2005sparse}. Sparse grid techniques have been incorporated in  collocation methods for   high-dimensional stochastic differential equations \cite{xiu2005high,xiu2007efficient,nobile2008sparse,ma2009adaptive},  Galerkin finite element methods \cite{zenger1991sparse, bungartz2004sparse,schwab2008sparse}, finite difference methods \cite{griebel1998adaptive, griebel1999adaptive}, finite volume methods \cite{hemker1995sparse}, and spectral methods \cite{griebel2007sparse, gradinaru2007fourier,shen2010sparse,shen2010efficient} for   high-dimensional PDEs.  However, their potential has not yet been fully realized under the DG framework, which constitutes the main focus of this paper.


Generally speaking, the design of  DG methods has two essential ingredients: (1) a (weak) formulation of the methods based on the underlying PDE, (2) a suitable approximation space. The appropriate choice of the formulation of the scheme as well as approximation properties of the discrete space can guarantee the methods with desirable properties such as stability and convergence.
 In fact, DG formulations based on traditional piecewise polynomial space have been extensively studied for the last few decades, and the methods are rather mature for various applications.
As for the   approximation spaces,   DG methods  are more flexible compared with continuous finite element methods  due to the lack of continuity requirement.
   In recent years, the ideas of using nonstandard spaces such as non-polynomial spaces \cite{yuan2006discontinuous, wang2009wkb}, or polynomial spaces with specific properties such as locally divergence-free properties \cite{Fengyan_Maxwell} have been explored. They are mainly driven by the needs to mimic  particular properties of the exact solution. In this paper, we aim at incorporating the sparse grid approximation spaces in the DG framework to treat high-dimensional problems.
   As an initial effort in our investigation,   we focus on high-dimensional elliptic problems and develop a sparse grid DG method for the following $d$-dimensional model problem
\begin{equation}
\label{elliptic eqn}
- \nabla \cdot (\bK  \nabla u) = f \  \ \ \text{in}\  \Omega=[0,1]^d,
\end{equation}
with a Dirichlet boundary condition
\begin{equation}
 u=g \  \ \ \text{on} \ \partial \Omega,
\end{equation}
where $u(\bx): \Omega \rightarrow \mathbb{R}$ is the unknown function and the  matrix-valued coefficient function $\bK = (k_{i,j})_{1\le i,j \le d}$  is symmetric positive definite and bounded below and above uniformly, i.e. there exist two positive constant $K_0, K_1$, such that
$$
\forall \,\bx \in \Omega, \quad K_0 \,\bx \cdot \bx \leq \bK \bx \cdot \bx \leq K_1 \,\bx \cdot \bx.
$$

Using a hierarchical basis representation, we construct a sparse finite element approximation space, reducing the degrees of freedom  from the standard {$O(h^{-d})$ to $O(h^{-1}|\log_2 h|^{d-1})$} for $d$-dimensional problems, where $h$ is the uniform mesh size in each dimension. Based on the IPDG framework and careful implementation, the new approximation space enables an efficient solver for
the model equation. \emph{A priori} error estimate shows that the proposed method is convergent in the order of $O(h^{k})$ up to the polylogarithmic term $|\log_2 h|^{d-1}$ in the energy norm when the solution is smooth enough, where $k$ is the degree of the polynomial space. Although the method  in this paper is restricted to problems  on a box-shaped domain, it has the potential to be implemented on general domains by considering either a coordinate transformation or a more complex sparse finite element space based on hierarchical decompositions on unstructured meshes.

The rest of this paper is organized as follows:
in Section \ref{sec:SGapproximation},  we  construct and analyze the DG approximation space on sparse grids.  In Section \ref{sec:sparseIP}, we formulate the scheme and perform error analysis of the sparse grid IPDG method. Implementation details will be discussed to ensure computational efficiency. In Section \ref{sec: SGDG numerical}, numerical examples in multi-dimensions (up to $d=4$) are provided to validate the accuracy and performance of the  method. Conclusions  and future work are given in Section \ref{sec:Con}.

\section{DG Finite Element Spaces on Sparse Grids}
\label{sec:SGapproximation}
In this section, we   introduce the main ingredient of our algorithm: the DG finite element space on sparse grids. We proceed in several steps. First, we  review the standard piecewise polynomial spaces in one dimension, and introduce a hierarchical decomposition with a set of orthonormal bases. Second, we   construct the sparse finite element space from  multi-dimensional multiwavelet bases as tensor products from the one-dimensional construction. Finally, we  discuss some key features and approximation properties of the sparse finite element space. Numerical tests are given to compare several definitions of sparse discontinuous finite element spaces.

\subsection{Hierarchical Decomposition of Piecewise Polynomial  Spaces in One Dimension}
\label{sec:1dmulti}
In this subsection, we will introduce the hierarchical representation of piecewise polynomial spaces in one dimension. Many discussions here are closely related to the original Haar wavelet \cite{haar}, the multiwavelet bases constructed in \cite{alpert1993class, alpert2002adaptive}, and multiresolution analysis (MRA) with DG methods for adaptivity
\cite{calle2005wavelets, archibald2011adaptive, iacono2011high, hovhannisyan2014adaptive, gerhard2013adaptive}
and constructions of trouble-cell indicator \cite{vuik2014multiwavelet} for hyperbolic conservation laws.

Without loss of generality, consider the interval $\Omega=[0,1]$, we define the $n$-th level grid $\Omega_n$, consisting of $I_{n}^j=(2^{-n}j, 2^{-n}(j+1)]$, $j=0, \ldots, 2^n-1$ with uniform cell size $h=2^{-n}.$ On this grid, we use $||\cdot||_{H^s(\Omega_n)}$ to denote the broken Sobolev norm, i.e. $||v||^2_{H^s(\Omega_n)}=\sum_{j=0}^{2^n-1} ||v||^2_{H^s(I_n^j)},$ where $||v||_{H^s(I_n^j)}$ is the standard Sobolev norm on $I_n^j.$ Similarly,  $|\cdot|_{H^s(\Omega_n)}$   denotes the broken Sobolev semi-norm, i.e. $|v|^2_{H^s(\Omega_n)}=\sum_{j=0}^{2^n-1} |v|^2_{H^s(I_n^j)}.$

We can define
$$V_n^k=\{v: v \in P^k(I_{n}^j),\, \forall \,j=0, \ldots, 2^n-1\}$$
to be the usual piecewise polynomials of degree at most $k$ on this grid.  Clearly, $V_n^k$ has degrees of freedom $2^n (k+1)$, and we notice that they have the nested structure for different values of $n,$
$$V_0^k \subset V_1^k \subset V_2^k \subset V_3^k \subset  \cdots$$
Due to this nested structure, we can define the multiwavelet subspace $W_n^k$, $n=1, 2, \ldots $ as the orthogonal complement of $V_{n-1}^k$ in $V_{n}^k$ with respect to the $L^2$ inner product on $\Omega$, i.e.
\begin{equation*}
V_{n-1}^k \oplus W_n^k=V_{n}^k, \quad W_n^k \perp V_{n-1}^k.
\end{equation*}

For notational convenience, we also denote the base space $W_0^k:=V_0^k$, which consists of all polynomials of up to degree $k$ on $[0,1]$. The dimension of $W_n^k$ is $2^{n-1}(k+1)$ when $n \geq 1$, and $k+1$ when $n=0$.
Now, we have found a hierarchical representation of  the standard piecewise polynomial space on grid $ \{I_n^j, j=0, \ldots, 2^n-1 \}$ as $V_N^k=\bigoplus_{0 \leq n\leq N} W_n^k$. We remark that   the space $W_n^k$, for $n>1$, represents the finer level details  when the mesh becomes refined and this is the key to the reduction in degrees of freedom in higher dimensions. 

For implementation purpose, we need to introduce basis functions for space $W_n^k$. The case of $n=0$ is trivial. It suffices to use a standard polynomial basis on $\Omega$. For example, by using the scaled Legendre polynomials, we can easily obtain a set of orthonormal bases in $W_0^k$.  When $n>1$, we will use the orthonormal bases introduced in \cite{alpert1993class}. In particular, the bases for $W_1^k$ are defined as
$$h_i(x)= 2^{{1}/{2}} f_i(2x-1), \quad i = 1, \ldots, k+1$$
where    $ \{f_i(x), \, i=1, \ldots, k+1 \}$ are functions supported on $[-1,1]$ depending on $k$, with the following properties.
\begin{enumerate}
\item The restriction of $f_i$ to the interval $(0, 1)$ is a polynomial of degree $k$.
\item The function $f_i$ is extended to $(-1, 0)$ as an even or odd function according to the parity of $i+k$:
$$f_i (x) = (-1)^{i+k} f_i (-x).$$
\item The bases have vanishing moments:
\begin{align*}
\int_{-1}^{1} f_j(x) x^i \,dx = 0 \quad i = 0, 1, \cdots, j+k-1.
 \end{align*}
\item The bases have the orthogonality and normality properties:
\begin{align*}
 \int_{-1}^{1} f_i(x) f_j(x) \,dx = <f_i, f_j> = \delta_{ij}, \quad i,j=1, \ldots, k.
 \end{align*}
\end{enumerate}

Functions $\{f_i\}$ can be computed by a repeated Gram-Schmidt algorithm. The particular form of $\{f_i\}$ up to $k=4$  are listed in Table 1 in \cite{alpert1993class}.
Multiwavelet bases in \cite{alpert1993class} retain the orthonormal properties of wavelet bases for different hierarchical levels by introducing the basis for $W_n^k, \,n \geq 1$  as

$$  v_{i,n}^j (x) = 2^{(n-1)/2} \,h_i(2^{n-1}x - j), \quad i = 1, \ldots, k+1,\, j=0, \ldots, 2^{n-1}-1.$$

To make the notations consistent and compact, when $n=0$, the bases $  v_{i,0}^0 (x), \, i=1, \ldots k+1$ are defined as the rescaled Legendre polynomials on $[0,1]$ with orthonormal property.
In summary, we have \cite{alpert1993class, alpert2002adaptive},
\begin{align*}
\int_{0}^{1}  v_{i,n}^j(x) v_{i',n'}^{j'}(x)\,dx = \delta_{i i'}\delta_{n n'}\delta_{j j'},
\end{align*}
where $\delta_{i i'}, \delta_{n n'}, \delta_{j j'}$ are the Kronecker delta symbols.

Next, we will discuss about the projection operators that are essential in   finite element analysis. We define $P_n^k$ as the standard univariate $L^2$  projection operator from  $L^2[0,1]$ to $V_n^k$. From there, we can introduce the
 increment projector,
\[Q_n^k: =
 \left\{
  \begin{array}{ll}
P_n^k - P_{n-1}^k,  & \textrm{if} \, n \ge1\\
P_0^k, & \textrm{if} \, n = 0,
  \end{array} \right.\]
  then we can see that
 $$W_n^k= Q_n^k L^2(0,1).$$

In summary, we arrive at the following identities for the hierarchical decomposition of the piecewise polynomial space and the projection operator:
  $$V_N^k=\bigoplus_{\substack{0 \leq n\leq N}} W_n^k,  \ P_N^k=\sum_{\substack{0 \leq n\leq N}} Q_n^k.$$

 The properties of $Q_n^k$ naturally rely on $P_n^k$.  For the $L^2$ projection operators, we recall the following approximation results. 

\begin{prop}[Convergence Property of the Projection Operator \cite{Ciarlet_1975_FEM_Elliptic}]
\label{prop2}

For a function $v \in  H^{p+1}(0,1)$,    we have the convergence property of the $L^2$ projection $P_n^k$ as follows: for any integer $t$ with $1 \leq t \leq \min\{p, k\}$,  and  $s\leq t+1$,
\begin{equation}
\label{eq:proj1}
 || P_n^k v- v||_{H^s(\Omega_n)}\le c_{k,s, t} 2^{-n(t+1-s)} ||v||_{H^{t+1}(\Omega)},
 \end{equation}
where $c_{k,s, t}$ is a constant that depends on $k, s, t$, but not on $n$. 
\end{prop}

From this property, using basic algebra, we can deduce that for $n \geq 1$,
$$ || Q_n^k v  ||_{H^s(\Omega_n)}\le \tilde{c}_{k,s, t} 2^{-n(t+1-s)} ||v||_{H^{t+1}(\Omega)}, $$
with
\begin{equation}
\label{eq:proj2}
\tilde{c}_{k,s, t}=c_{k,s, t} \,(1+2^{-(t+1-s)}).
 \end{equation}

  Finally, we want to remark on a subtle issue concerning the hierarchical decomposition and the multiwavelet space. If we use a projection other than $P_n^k$, say $\underline{P}_n^k$ is another projection from $L^2(\Omega)$ to $V_n^k.$ Then the increment projector $\underline{Q}_n^k \ne Q_n^k$. The space $\underline{W}_n^k=\underline{Q}_n^k L^2(0,1)$ is different from $W_n^k,$ and likewise for the multiwavelet bases. However, the hierarchical structure $V_{n-1}^k \oplus \underline{W}_n^k=V_{n}^k$ still holds. In our discussion in the next subsection, we will highlight the implication of this statement and show that using different definitions of the projector and the increment space will not affect the definition of  the sparse discontinuous finite element space in multi-dimensions.


\subsection{Sparse Discontinuous Finite Element Spaces in Multi-dimensions}
\label{sec:multibases}

In this subsection,  we   introduce the sparse finite element space constructed from  multi-dimensional multiwavelet bases obtained from tensor products of the one-dimensional bases in Section \ref{sec:1dmulti}.

For a $d$-dimensional problem, we consider the domain $\bx=(x_1, \ldots, x_d) \in \Omega=[0,1]^d$.   To facilitate the discussion, we first introduce some notations on the norms and operations of multi-indices in $\mathbb{N}_0^d$, where $\mathbb{N}_0$ denotes the set of nonnegative integers. For $\ba =(\alpha_1, \ldots, \alpha_d) \in \mathbb{N}_0^d,$ we define the $l^1$ and $l^\infty$ norms
$$
|\ba|_1:=\sum_{m=1}^d \alpha_m, \qquad   |\ba|_\infty:=\max_{1\leq m \leq d} \alpha_m,
$$
the component-wise arithmetic operations
$$
\ba \cdot \bb :=(\alpha_1 \beta_1, \ldots, \alpha_d \beta_d), \qquad c \cdot \ba:=(c \alpha_1, \ldots, c \alpha_d), \qquad 2^\ba:=(2^{\alpha_1}, \ldots, 2^{\alpha_d}),
$$
and the relational operators
$$
\ba \leq \bb \Leftrightarrow \alpha_m \leq \beta_m, \, \forall m
$$
$$
\ba<\bb \Leftrightarrow \ba \leq \bb \textrm{  and  } \ba \neq \bb.
$$

Now, for a multi-index $\bl=(l_1, \ldots, l_d) \in \mathbb{N}_0^d$, which indicates the level of the mesh in a multivariate sense, we consider the standard rectangular grid $\Omega_\bl$ with mesh size
$h_\bl:=(2^{-l_1}, \ldots, 2^{-l_d}).$ We define the smallest size among all dimensions to be $h=\min\{2^{-l_1}, \ldots, 2^{-l_d} \}$,
  an elementary cell $I_\bl^\bj=\{\bx: x_m \in (2^{-l_m}j_m, 2^{-l_m}(j_m+1)\}$,   and
$$\bV_\bl^k=\{\bv: \bv(\bx) \in P^k(I^{\bj}_{\bl}), \,\,  \bzero \leq \bj  \leq 2^{\bl}-\bone \},  $$
where $P^k(I^{\bj}_{\bl})$ denotes polynomials of degree up to $k$ in each dimension on cell $I^{\bj}_{\bl}$. The space $\bV_\bl^k$ contains the traditional tensor-product piecewise polynomials used in the DG discretizations. Moreover, if we use equal refinement of size $2^{-N}$ in each direction, we denote the space to be $\bV_N^k$, and it consists of $(2^{N}(k+1))^d$ degrees of freedom.


The foundation of sparse grids is to use the tensor product of  the one-dimensional hierarchical decomposition, and only chooses the most relevant bases to guarantee suitable approximation properties. To illustrate the main ideas, we introduce

$$\bW_\bl^k=W_{l_1, x_1}^k \times W_{l_2, x_2}^k \cdots \times W_{l_d, x_d}^k,$$
where $W_{l_m, x_m}^k$ corresponds to the space $W_{l_m}^k$ defined in the $m$-th dimension as defined in the previous subsection.  The number of degrees of freedom associated with $\bW_\bl^k$ is
$$dim(\bW_\bl^k )= \prod_{m=1}^d dim(W_{l_m}^k).$$
Based on the one-dimensional hierarchical decomposition, it is easy to see
$$\bV_\bl^k=V_{l_1,x_1}^k \times V_{l_2,x_2}^k \cdots \times V_{l_d,x_d}^k=\bigoplus_{j_1 \leq l_1, \ldots, j_d \leq l_d} \bW_\bj^k,$$
and
$$\bV_N^k=V_{N,x_1}^k \times V_{N,x_2}^k \cdots \times V_{N,x_d}^k=\bigoplus_{|\bl|_\infty \leq N} \bW_\bl^k.$$

The basis functions for $\bW_\bl^k$ can be defined by a tensor product
$$v_{\bi, \bl}^{\bj}(\bx) := \prod_{m=1}^d {v_{i_m, l_m}^{j_m}(x_m)}, \quad j _m= 0,\ldots, \max(0, 2^{l_m-1}-1);\,  i_m= 1, \ldots, k+1.$$
They form a set of orthonormal basis due to the property of the one-dimensional bases.

 Now, we are ready to introduce the sparse finite element approximation space. In particular, we define
$$ \hat{\bV}_N^k:=\bigoplus_{\substack{ |\bl|_1 \leq N} }\bW_\bl^k.$$
This definition is motivated by continuous finite element space on sparse grid \cite{bungartz2004sparse}. Instead of considering the standard piecewise polynomial space $\bV_N^k,$ which has exponential dependence on the dimension $d$, we shall use its subspace  $\hat{\bV}_N^k.$ This space has good approximation results (see Section \ref{sec:sparseapprox} for details) with significantly reduced degrees of freedom. We remark that the key in the construction lies in the choice of condition for $\bl$. Here we have taken it to be $|\bl|_1 \leq N$. As demonstrated in \cite{bungartz2004sparse}, other conditions are also possible to realize dimension reduction with aims of optimizing various types of approximation results. We defer the numerical study for the comparison of several different spaces to Section \ref{sec:approx}.

\begin{rem}
We mentioned in the previous subsection that there are different ways to construct the increment space $W_n^k.$ However, the choice  will not affect the definition of $\hat{\bV}_N^k.$ This is because
$$ \bigoplus_{\substack{ |\bl|_1 \leq N} }\bW_\bl^k=\bigoplus_{\substack{ |\bl|_1 \leq N} }\bV_\bl^k=\bigoplus_{\substack{ |\bl|_1 \leq N} }\underline{\bW}_\bl^k,$$
for any space $\underline{\bW}_\bl^k$ constructed from one dimensional increment space $\underline{W}_n^k$ satisfying $V_{n-1}^k \oplus \underline{W}_n^k=V_{n}^k.$
\end{rem}




The next lemma will give a count of dimensions for the space $ \hat{\bV}_N^k.$
\begin{lem}
\label{lem:count}
The dimension of $ \hat{\bV}_N^k$ is given by
\begin{flalign*}
 &dim(\hat{\bV}_N^k)= (k+1)^d \left \{\sum_{m=0}^{d-1} \binom{d}{m} \left((-1)^{d-m}+ 2^{N+m-d+1}  \sum_{n=0}^{d-m-1} \binom{N}{n} \cdot (-2)^{d-m-1-n} \right)+1 \right \}.
\end{flalign*}
 Suppose there is an upper bound on the dimension $d \leq d_0$, then there exist  constants $c_{d_0}, C_{d_0}$ depending only on $d_0$, such that
\begin{flalign*}
 &c_{d_0}  (k+1)^d 2^N N^{d-1} \leq dim(\hat{\bV}_N^k)\leq  C_{d_0}  (k+1)^d 2^N N^{d-1}.
\end{flalign*}
 \end{lem}

\noindent{\it Proof}: Due to the distinctiveness of the zero-th level in a mesh, we need to  distinguish the zero-th level from other levels in the proof.
  Therefore, for each multi-index $\bl$, we define  the set of the dimensions with  mesh level 0 as $\bl^0 = \{ i| \,l_i = 0,  i=1, \ldots, d \} $
 and the number of such dimensions as  $$|\bl |_0:= \# \textrm{ in } \bl^0.$$
 Then
\begin{flalign*}
 dim(\hat{\bV}_N^k)&=dim(\bigoplus_{\substack{ |\bl|_1 \leq N} } \bW_\bl^k) \notag \\
 &=\sum_{m=0}^d dim(\bigoplus_{\substack{ |\bl|_0=m, |\bl|_1 \leq N }} \bW_\bl^k).
\end{flalign*}

 We will discuss two cases based on the value of $|\bl|_0.$

 {\bf Case 1. $0 \le |\bl|_0=m \le d-1$}, i.e., there are $m$ dimensions of level zero from the $d$ dimensions of the multi-index $\bl$. Clearly,

\begin{flalign*}
dim(  \bigoplus_{\substack{  |\bl|_0=m, |\bl|_1 \leq N}} \bW_\bl^k )
= \binom {d}{m} \,\,dim(  \bigoplus_{\substack{ l_t = 0, \, 1 \le t \le m; \, l_t> 0, \, t >m\\|\bl|_1 \leq N }} \bW_\bl^k ).
\end{flalign*}

Since there are $d-m$ dimensions that have meshes of level no less than one, we always have $  |\bl |_1 \ge d-m$.
Define $e_m = (0,\cdots, 0, 1, \cdots, 1)$ to be a vector in  $\mathbb{N}_0^d$ whose first $m$ dimensions are $0$, and the rest $d-m$ dimensions are $1$. Then, by Lemma 3.6 in \cite{bungartz2004sparse},
\begin{flalign*}
dim( \bigoplus_{\substack{ l_t = 0, \, 1 \le t \le m; \, l_t> 0, \, t >m \\ |\bl|_1 \leq N }} \bW_\bl^k )
 =& (k+1)^d\sum_{\substack{ l_t = 0, \, 1 \le t \le m; \, l_t> 0, \, t >m \\ |\bl|_1 \leq N }} {2^{|\bl - {e_m}|_1} }\\
 = & (k+1)^d\sum_{n=d-m}^{N}{2^{n-(d-m)} \sum_{\substack{ l_t = 0, \, 1 \le t \le m; \, l_t> 0, \, t >m \\ |\bl|_1 =n }} 1}\\
= & (k+1)^d\sum_{n=d-m}^{N}{2^{n-(d-m)} \binom {n-1}{(d-m)-1} }\\
= &(k+1)^d\sum_{n=0}^{N+m-d}{2^{n} \binom {n+(d-m)-1}{(d-m)-1} } \\
= & (k+1)^d \left( (-1)^{d-m}+ 2^{N+m-d+1}   \sum_{n=0}^{d-m-1} \binom{N}{n} \cdot (-2)^{d-m-1-n} \right).
\end{flalign*}

{\bf Case 2. $|\bl|_0= d$} means $\bl=\bzero.$  Therefore,
\begin{flalign*}
 dim(\bigoplus_{\substack{ |\bl|_0=d, |\bl|_1 \leq N }} \bW_\bl^k)= dim(  \bW_\bzero^k)=(k+1)^d.
\end{flalign*}

Finally,  we combine both cases, and arrive at
\begin{flalign*}
 &dim(\hat{\bV}_N^k)=\sum_{m=0}^d dim(\bigoplus_{\substack{ |\bl|_0=m, |\bl|_1 \leq N }} \bW_\bl^k) \\
 &= (k+1)^d \left \{\sum_{m=0}^{d-1} \binom{d}{m} \left((-1)^{d-m}+ 2^{N+m-d+1}  \sum_{n=0}^{d-m-1} \binom{N}{n} \cdot (-2)^{d-m-1-n} \right)+1 \right \}.
\end{flalign*}

We notice that
\begin{flalign*}
&\sum_{m=0}^{d-1} \binom{d}{m}  2^{N+m-d+1} \sum_{n=0}^{d-m-1} \binom{N}{n} \cdot (-2)^{d-m-1-n} \\
& \leq \sum_{m=0}^{d-1} \binom{d}{m}  2^{N+m-d+1} C_{d_0} N^{d-m-1}\\
& \leq C_{d_0} 2^{N-d+1} N^{d-1}  \sum_{m=0}^{d} \binom{d}{m}  2^{ m}  N^{-m}\\
&= C_{d_0} 2^N N^{d-1} (1+\frac{2}{N})^d \leq C_{d_0} 2^N N^{d-1},
\end{flalign*}
where we use $C_{d_0}, c_{d_0}$ to denote  generic constants that depend only on $d_0$, not on $d, k, N$. They may have different values in each occurrence throughout the proof.

Similarly,
\begin{flalign*}
 &dim(\hat{\bV}_N^k) \ge dim(\bigoplus_{\substack{ |\bl|_0=0, |\bl|_1 \leq N }} \bW_\bl^k)=(k+1)^d \left( (-1)^{d}+ 2^{N-d+1}   \sum_{n=0}^{d-1} \binom{N}{n} \cdot (-2)^{d-1-n} \right) \\
 &\ge (k+1)^d (-1+c_{d_0}2^{N-d+1}N^{d-1}) \geq c_{d_0} (k+1)^d  2^N N^{d-1}.
\end{flalign*}

In summary, there exist  constants $c_{d_0}, C_{d_0}$, such that
\begin{flalign*}
 &c_{d_0}  (k+1)^d 2^N N^{d-1} \leq dim(\hat{\bV}_N^k)\leq  C_{d_0}  (k+1)^d 2^N N^{d-1},
\end{flalign*}
and we are done.
 \hfill $ \blacksquare $

 This lemma implies that, upon mesh refinement, the degrees of freedom for the sparse finite element space will grow in the order of $h^{-1} |\log_2{h}|^{d-1}$ instead of $h^{-d}$ for the traditional piecewise polynomial space. This translates into a significant reduction in computational cost when $d$ is large.  The next subsection will verify approximation property of the space  $\hat{\bV}_N^k.$

%

%

\subsection{Approximation Results}
\label{sec:sparseapprox}


Given any one-dimensional projection operator  $\underline{P}^k_n$ and its increment operator $\underline{Q}^k_n$, we can naturally define a projection from $H^1(\Omega)$ to $\hat{\bV}_N^k$ as,
\begin{equation}
\label{eqn:multiproj}
\underline{\hat{\bP}}_N^k:=\sum_{\substack{ |\bl|_1 \leq N}} \underline{Q}^k_{l_1, x_1} \otimes \cdots \otimes \underline{Q}^k_{l_d, x_d},
\end{equation}
where $\underline{Q}^k_{l_m, x_m}$ denotes the operator $\underline{Q}^k_{l_m}$ in the $m$-th dimension. In \cite{schwab2008sparse}, continuous sparse finite element space is used with streamline diffusion method to solve transport-dominated diffusion problem. The finite element space considered is $\hat{\bV}_N^{c,k}:=\hat{\bV}_N^k \cap C^0(\Omega)$ for $k \geq 1.$ The approximation properties for   projector defined in \eqref{eqn:multiproj} in $L^2(\Omega), H^1(\Omega), H^2(\Omega_N)$ norms are obtained when   $\underline{P}^k_n$ is a  univariate  projector onto $V_n^k \cap C^0([0,1])$
satisfying \eqref{eq:proj1}. Because $\hat{\bV}_N^{c,k}$ is a subset of $\hat{\bV}_N^k,$ we can directly use
 the results in \cite{schwab2008sparse} and obtain an estimate of the projection error for $\underline{\hat{\bP}}_N^k$.

We introduce some notations about the seminorm of a function's mixed derivatives.
For any set $I=\{i_1, \ldots i_r \} \subset \{1, \ldots d\}$ with $|I|$=r, we define $I^c$ to be the complement set of $I$ in $\{1, \ldots d\}$. For non-negative integers $\alpha, \beta$ and set $I$,  we define the seminorm
\begin{flalign*}
&|v|_{H^{\alpha, \beta, I}(\Omega)} := \\
&\sum_{0 \le \alpha_1 \le \alpha} \cdots \sum_{0 \le \alpha_r \le \alpha} \sum_{0 \le \beta_1 \le \beta} \cdots \sum_{0 \le \beta_{d-k} \le \beta} \left | \left |\left ( \frac{\partial^{\alpha_1}}{\partial x_{i_1}^{\alpha_1}} \cdots \frac{\partial^{\alpha_r}}{\partial x_{i_r}^{\alpha_r}}  \right ) \left ( \frac{\partial^{\beta_1}}{\partial x_{j_1}^{\beta_1}} \cdots \frac{\partial^{\beta_{d-r}}}{\partial x_{j_{d-r}}^{\beta_{d-r}}}  \right )v \right | \right|_{L^2(\Omega)},
\end{flalign*}
and 
$$
|v|_{\mathcal{H}^{t+1}(\Omega)} :=\max_{s \in \{0, 1\}} \max_{1 \leq r \leq d} \left ( \max_{|I|=r} |v|_{H^{t+1, s, I}(\Omega)} \right ).
$$
We recall the following results in \cite{schwab2008sparse}  about the approximation properties   in $L^2(\Omega)$, $H^1(\Omega)$ and  $H^2(\Omega_N)$ norm.



\begin{thm}
\label{thm:appx}
Let $\underline{\hat{\bP}}_N^k$  defined in \eqref{eqn:multiproj} be constructed from $\underline{P}^k_n$ which is a  univariate  projector onto $V_n^k \cap C^0([0,1])$
satisfying \eqref{eq:proj1}, then for  $k \geq 1$, any $1 \leq t \leq \min \{p, k\}$, there exist constant $\underline{c}_{k,t}, \,\kappa_0(k, t, N), \,\kappa_1(k, t, N) >0$, such that for any
$v \in \mathcal{H}^{p+1}(\Omega)$, $N\geq 1$, $d \geq 2$,  we have
$$
  | \underline{\hat{\bP}}_N^k v- v |_{H^s(\Omega_N)}=| \underline{\hat{\bP}}_N^k v- v |_{H^s(\Omega)}\le  \underline{c}_{k,t} d^{1+s/2} \kappa_s(k, t,  N)^{d-1+s} 2^{-N(t+1-s)} |v|_{\mathcal{H}^{t+1}(\Omega)},
$$
for $s=0, 1$ ($L^2$ norm and $H^1$ seminorm), and
  \[  | \underline{\hat{\bP}}_N^k v- v |_{H^2(\Omega_N)} \leq
 \left\{
  \begin{array}{ll}
\underline{c}_{k,t} \left (d^{3/2}\kappa_1(k, t, N)^d +d^2 \kappa_0(k, t, N)^{d-1} 2^{-N}\right) 2^{-N(t-1)} |v|_{\mathcal{H}^{t+1}(\Omega)},  & k \ge 2\\
\underline{c}_{k,t} \left (d^{1/2}  +d^2 \kappa_0(k, t,  N)^{d-1} 2^{-N} \right)  |v|_{\mathcal{H}^{2}(\Omega)}, & k=1,
  \end{array} \right.\]
where
\[\kappa_s(k, t, N) =
 \left\{
  \begin{array}{ll}
\tilde{c}_{k,0,t}(N+1)e^{1/(N+1)}+c_{k,0,0},  & s=0\\
2\tilde{c}_{k,0,t}+c_{k,0,0}, & s=1,
  \end{array} \right.\]
  and  $\tilde{c}_{k,0,t}, c_{k,0,0}$  are defined in \eqref{eq:proj1} and \eqref{eq:proj2}.
Moreover, for $s=0$, if
  $$
  \tilde{c}_{k,0,t} 2^{t+1}/(2^{t+1}-1)+c_{k,0,0}<1,
  $$
  there exists a positive constant $c_{t,k}$, such that $\kappa(k, t, s, N)<1$ for all $N \geq 1, d \geq 2$ with $N+1 \leq c_{t,k}(d-1).$
\end{thm}

This theorem implies a convergence rate of $O(h^{k+1})$ up to the polylogarithmic term $|\log_2 h|^{d-1}$ in the $L^2$ norm, which is comparable with the traditional full grid approach when the function $v$ retains enough smoothness (we actually do require higher regularity of the function $v$ when compared with standard piecewise polynomial spaces). This estimate  serves as the foundation in the error estimates of IPDG method in the energy norm. We remark that when switching the projection operator (say $L^2$ projection or other projectors on to $\hat{\bV}_N^k$) the error bound may change. However, our estimates in Section \ref{sec:error_analysis} do not require any particular property from the projection. Therefore, it suffices for us to use this theorem directly from \cite{schwab2008sparse}. In future work, we will establish convergence property of different projections onto $\hat{\bV}_N^k,$ for $k \geq 0$. This will be   essential for error estimates for broader definitions of DG methods for   other types of equations.



%
%
%

\subsection{Approximations from Various  Sparse Discontinuous Finite Element Spaces: A Numerical Investigation}
\label{sec:approx}

In this subsection, we   investigate the approximation properties of three sparse finite element spaces from a numerical perspective. As one can imagine, there is an intricate relation between the number of  bases used with the approximation properties of the spaces. Indeed, the choice adopted in Section \ref{sec:multibases} is not unique in sparse grids literature. Ideally, how to choose the ``best" bases is  an optimization problem of cost/benefit ratios  \cite{griebel2005sparse}, and the ultimate strategy may be an adaptive algorithm. Our preliminary study in this subsection compares three   intuitive choices of spaces,
 including
$\hat{\mathbf{V}}_N^k$ defined in the previous subsections, and
$$ \tilde{\bV}_N^k:=\bigoplus_{\substack{ |\bl|_1 \leq N+d-1 \\ |\bl|_\infty\leq N }}\bW_\bl^k,$$
which is a slightly larger space. $\tilde{\bV}_N^k$ contains   the sparse continuous piecewise linear space in \cite{bungartz2004sparse} when $k \geq 1$. Moreover,   by a similar argument as in the proof of Lemma \ref{lem:count},    its degrees of freedom, though slightly larger, are also of the order $h^{-1} |\log_2{h}|^{d-1}$. 
%

The other space we consider retains not only sparsity with respect to $h$, but also with respect   to the degree of polynomials (we call this  $p$-sparsity).
If we introduce
$$\hat{\bW}_\bl^k=\bigoplus_{|\bk|_1\leq k}W_{l_1, x_1}^{k_1} \times W_{l_2, x_2}^{k_2} \cdots \times W_{l_d, x_d}^{k_d},$$
with $\bk=(k_1,\cdots,k_d)$,  the new sparse space can be defined as
 $$\hat{\hat{\mathbf{V}}}_N^k:=\bigoplus_{\substack{ |\bl|_1 \leq N} }\hat{\bW}_\bl^k.$$
 It is clear that
 $$\hat{\hat{\mathbf{V}}}_N^k\subset\hat{\mathbf{V}}_N^k\subset\tilde{\bV}_N^k.$$

A similar argument of  Lemma \ref{lem:count}   gives,
\begin{lem}
The dimension of $ \hat{\hat{\bV}}_N^k$ is
\begin{flalign*}
 &dim(\hat{\hat{\bV}}_N^k)= \binom{k+d}{d}
  \left \{\sum_{m=0}^{d-1} \binom{d}{m} \left((-1)^{d-m}+ 2^{N+m-d+1}  \sum_{n=0}^{d-m-1} \binom{N}{n} \cdot (-2)^{d-m-1-n} \right)+1 \right \}.
\end{flalign*}
 Suppose there is an upper bound on the dimension $d \leq d_0$, then there exist  constants $c_{d_0}, C_{d_0}$ depending only on $d_0$, such that
\begin{flalign*}
 &c_{d_0}  \binom{k+d}{d} 2^N N^{d-1} \leq dim(\hat{\hat{\bV}}_N^k)\leq  C_{d_0}   \binom{k+d}{d} 2^N N^{d-1}.
\end{flalign*}
 \end{lem}

Note that, when $k$ and $d$ are large, $dim(\hat{\hat{\bV}}_N^k)$ can be significantly smaller than $dim(\hat{\bV}_N^k)$ and $dim(\tilde{\bV}_N^k)$.
 As for the implementation of the space $ \hat{\hat{\bV}}_N^k,$ we can no longer use the orthogonal hierarchical bases  in Section \ref{sec:1dmulti} to represent functions in $\hat{\hat{\mathbf{V}}}_N^k$, since the vanishing-moment functions $f_i$, which are used to construct the orthogonal hierarchical bases, have the same degrees.
  Instead, we can adopt another set of non-orthogonal hierarchical basis functions defined as follows. We first construct the bases for $W_n^k$.
 Again, the case of $n=0$ is trivial. Similar to the orthogonal case, the bases for $W_1^k$ are defined as
 $$\hat{h}_i(x) = 2^{1/2}\hat{f}_i(2x-1),\quad i=1,\cdots,k+1,$$
 where $\hat{f}_i,\,i=1,\cdots,k+1$ are functions supported on $[-1,1]$ and defined by,
  $$\hat{f}_i(x) = x^{i-1},\quad x\in(-1,0)$$ and
 $$\hat{f}_i(x) = -x^{i-1},\quad x\in(0,1).$$
 The bases for $W_n^k,\,n\ge1$ are
 $$  \hat{v}_{i,n}^j (x) = 2^{(n-1)/2} \,\hat{h}_i(2^{n-1}x - j), \quad i = 1, \ldots, k+1,\, j=0, \ldots, 2^{n-1}-1.$$
In summary, the basis functions in $\hat{\bW}_\bl^k$ (and $\hat{\hat{\bV}}_N^k$) are
 $$\hat{v}_{\bi, \bl}^{\bj}(\bx) := \prod_{\substack{m=1\\ |\bi|_1\leq k}}^d {v_{i_m, l_m}^{j_m}(x_m)}, \quad j _m= 0,\ldots, \max(0, 2^{l_m-1}-1).$$
This set of basis is no longer orthonormal in the $L^2$ sense, but a Gram-Schmidt algorithm can be adopted to achieve the purpose.

Ideally, we should be able to use detailed analysis to compare the approximation results from the three spaces. However, for simplicity of implementation, here we only   provide a numerical comparison for the approximation of a smooth function
$$u(\mathbf{x})=\exp\left(\prod_{m=1}^d x_i\right),\quad \mathbf{x}\in[0,1]^d,$$
by the standard $L^2$ projections:
$\hat{\bP}^k_n$, $\tilde{\bP}^k_n$, and $\hat{\hat{\bP}}^k_n$   onto $\hat{\mathbf{V}}^k_N$, $\tilde{\mathbf{V}}^k_N$, and $\hat{\hat{\mathbf{V}}}_N^k$, respectively. We measure various norms of the projection errors and the degrees of freedom for each space to find the best balance between accuracy and efficiency among the three choices.

In Tables \ref{table:approx_d2p2}-\ref{table:approx_d3p2}, we report the projection errors $e_1=\hat{\bP}^k_nu-u$, $e_2=\tilde{\bP}^k_nu-u$ and $e_3=\hat{\hat{\bP}}^k_nu-u$ in $L^1$, $L^2$, $L^\infty$, and $H^1$ norms and the associated orders of accuracy for $k=2$, $d=2, 3$. The detailed computational procedures as well as errors for other values of $d$ and $k$ are reported in \cite{zixuanthesis}.
%
For all three spaces, the degrees of freedom are all significantly less than $(k+1)^d 2^{Nd}$, which is the degree of freedom for the full grid approximation. For  space $ \tilde{\bV}_N^k$, $(k+1)$-th order of accuracy is clearly observed for the $L^1$ and $L^2$ errors, and $k$-th order of accuracy is observed for the $H^1$ error, while slight reduction of accuracy is observed for the $L^\infty$ error. For space $ \hat{\bV}_N^k$, $k$-th order of accuracy is observed for the $H^1$ error and slight reduction of accuracy is observed for $L^1$ and $L^2$ errors. However, we do observe about half-order to one order reduction of accuracy for the $L^\infty$ error. We remark that  the magnitude of errors computed by space $ \hat{\bV}_N^k$ is larger than that by space $ \tilde{\bV}_N^k$ for a fixed $N$. The apparent reason is that the degrees of freedom of $ \hat{\bV}_N^k$ are smaller than $ \tilde{\bV}_N^k$. However, when the degrees of freedom are comparable for the two spaces, the magnitude of errors computed by space  $ \hat{\bV}_N^k$ is smaller than $ \tilde{\bV}_N^k$. Therefore, $ \hat{\bV}_N^k$ is preferred if the computational efficiency is concerned. 
For $\hat{\hat{\mathbf{V}}}_N^k$, which has the least degrees of freedom among the three sparse spaces, severe reduction of accuracy is observed for all error norms. Moreover, the magnitude of errors is significantly larger than that by space  $ \hat{\bV}_N^k$ or $ \tilde{\bV}_N^k$ with comparable degrees of freedom.

 In summary, based on the discussion, we will adopt space $\hat{\mathbf{V}}_N^k$ in the computation for elliptic equations in  the next section. A further reduction in computational cost is possible by finding the optimal subset of $\hat{\mathbf{V}}_N^k$ by adaptive algorithms. We do not pursue this direction in this paper, and leave it to future study.

\begin{table}
\caption{Projection errors and orders of accuracy. SGDOF denotes the degrees of freedom of the   sparse approximation space in use.
FGDOF denotes the degrees of freedom of the full grid approximation space and is equal to $2^{Nd}(k+1)^d.$
$d=2$. $k=2$.}
\vspace{2 mm}
\centering
\begin{tabular}{c c c c c c c c c c c c}
\hline

N & SGDOF & FGDOF  & $L^1$ error & order & $L^2$ error & order & $L^\infty$ error & order  &$H^1$ error & order \\
\hline
     &  & &  &  &$\tilde{\mathbf{V}}^k_N$  &   & &  &  & \\
\hline
2	&	108	&	144 &2.82E-05	&		&	4.36E-05	&		&	4.05E-04	&		&	2.26E-03	&		\\
3	&	432	&	576& 3.50E-06	&	3.01	&	5.47E-06	&	2.99	&	5.94E-05	&	2.77	&	5.66E-04	&	2.00	\\
4	&	720	&	2304& 4.39E-07	&	3.00	&	6.86E-07	&	2.99	&	8.54E-06	&	2.80	&	1.42E-04	&	2.00	\\
5	&	1728	&9216 &	5.51E-08	&	3.00	&	8.61E-08	&	3.00	&	1.21E-06	&	2.81	&	3.54E-05	&	2.00	\\
6	&	4032	&36864 &	6.90E-09	&	3.00	&	1.08E-08	&	3.00	&	1.71E-07	&	2.83	&	8.84E-06	&	2.00	\\

\hline
     &  & &  & & $\hat{\mathbf{V}}^k_N$   &  & &  &  & \\
\hline
2	&	72	&	144 &3.51E-05	&		&	5.23E-05	&		&	6.48E-04	&		&	2.41E-03	&		\\
3	&	180	&	576&4.84E-06	&	2.86	&	7.26E-06	&	2.85	&	1.23E-04	&	2.40	&	6.08E-04	&	1.99	\\
4	&	432	&	2304&6.56E-07	&	2.88	&	9.96E-07	&	2.87	&	2.21E-05	&	2.47	&	1.53E-04	&	1.99	\\
5	&	1008	&	9216 &8.81E-08	&	2.90	&	1.35E-07	&	2.88	&	3.77E-06	&	2.55	&	3.82E-05	&	2.00	\\
6	&	2304	&36864 &	1.17E-08	&	2.91	&	1.81E-08	&	2.90	&	6.13E-07	&	2.62	&	9.55E-06	&	2.00	\\

\hline

     &  & & & & $\hat{\hat{\mathbf{V}}}^k_N$    &  & &  &  & \\
\hline
2	&	48	&144 &	1.73E-03	&		&	2.51E-03	&		&	2.78E-02	&		&	5.62E-02	&		\\
3	&	120	&576&	4.86E-04	&	1.83	&	7.22E-04	&	1.80	&	1.09E-02	&	1.35	&	2.80E-02	&	1.01	\\
4	&	288	&2304&	1.35E-04	&	1.85	&	2.02E-04	&	1.84	&	3.99E-03	&	1.46	&	1.39E-02	&	1.01	\\
5	&	672	&	9216 &	3.68E-05	&	1.88	&	5.53E-05	&	1.87	&	1.37E-03	&	1.55	&	6.95E-03	&	1.00	\\
6	&	1536	&36864 &	9.90E-06	&	1.89	&	1.49E-05	&	1.89	&	4.45E-04	&	1.62	&	3.47E-03	&	1.00	\\
\hline

\end{tabular}
\label{table:approx_d2p2}

\end{table}

\begin{table}
\caption{Projection errors and orders of accuracy. SGDOF denotes the degrees of freedom of the sparse approximation space in use.
FGDOF denotes the degrees of freedom of the full grid approximation space  and is equal to $2^{Nd}(k+1)^d.$
$d=3$. $k=2$.}
\vspace{2 mm}
\centering
\begin{tabular}{c c c c c c c c c c c c}
\hline

N & SGDOF& FGDOF  & $L^1$ error & order & $L^2$ error & order & $L^\infty$ error & order  &$H^1$ error & order \\
\hline
     &  & &  &  & $\tilde{\mathbf{V}}^k_N$ &    & &  &  & \\
\hline
2	&	1188	&1728&	8.65E-06	&		&	1.88E-05	&		&	5.64E-04	&		&	9.83E-04	&		\\
3	&	4104	&13824&	1.08E-06	&	3.00	&	2.36E-06	&	3.00	&	8.11E-05	&	2.80	&	2.45E-04	&	2.00	\\
4	&	12096	&110592&	1.35E-07	&	3.00	&	2.95E-07	&	3.00	&	1.15E-05	&	2.82	&	6.12E-05	&	2.00	\\
5	&	32832	&884736&	1.69E-08	&	3.00	&	3.69E-08	&	3.00	&	1.65E-06	&	2.80	&	1.53E-05	&	2.00	\\
6	&	84672	&7077888&	2.12E-09	&	3.00	&	4.62E-09	&	3.00	&	2.40E-07	&	2.78	&	3.83E-06	&	2.00	\\

\hline

     &  & & & & $\hat{\mathbf{V}}^k_N$  &   & &  &  & \\
\hline
2	&	351	&	1728&	1.41E-05	&		&	2.58E-05	&		&	1.33E-03	&		&	1.10E-03	&		\\
3	&	1026	&	13824&2.12E-06	&	2.73	&	3.86E-06	&	2.74	&	3.16E-04	&	2.08	&	2.80E-04	&	1.98	\\
4	&	2808	&110592&	3.15E-07	&	2.75	&	5.76E-07	&	2.74	&	7.07E-05	&	2.16	&	7.07E-05	&	1.99	\\
5	&	7344	&884736&	4.62E-08	&	2.77	&	8.56E-08	&	2.75	&	1.50E-05	&	2.24	&	1.77E-05	&	1.99	\\
6	&	18576	&7077888&	6.66E-09	&	2.79	&	1.26E-08	&	2.76	&	3.01E-06	&	2.31	&	4.44E-06	&	2.00    \\

\hline

     &  & & & & $\hat{\hat{\mathbf{V}}}^k_N$   &  & &  &  & \\
\hline

2	&	130&	1728	&	4.50E-03	&		&	6.59E-03	&		&	1.13E-01	&		&	9.32E-02	&		\\
3	&	380	&	13824&	1.82E-03	&	1.30	&	2.71E-03	&	1.28	&	7.36E-02	&	0.62	&	5.93E-02	&	0.65	\\
4	&	1040&110592	&	7.51E-04	&	1.28	&	1.11E-03	&	1.29	&	4.30E-02	&	0.78	&	3.91E-02	&	0.60	\\
5	&	2720&884736	&	3.10E-04	&	1.28	&	4.53E-04	&	1.29	&	2.32E-02	&	0.89	&	2.70E-02	&	0.54	\\
6	&	6880&7077888	&	1.30E-04	&	1.25	&	1.88E-04	&	1.26	&	1.18E-02	&	0.97	&	1.95E-02	&	0.47	\\

\hline

\end{tabular}
\label{table:approx_d3p2}

\end{table}

\section{IPDG Method on Sparse Grids for Elliptic Equations}
\label{sec:sparseIP}

In this section, we solve  the   second order linear elliptic boundary value problem \eqref{elliptic eqn} by IPDG methods with the sparse finite element space $\hat{\bV}^k_N.$ We will first formulate the scheme and discuss its implementation details, and then perform \emph{a priori} error analysis.

\subsection{Formulation of the Scheme}
\label{sec:formulation}
First, we   introduce some basic notations about jumps and averages for   piecewise functions defined on a grid $\Omega_N$. Let $\Gamma:=\bigcup_{T \in \Omega_N} \partial_T$ be the union of the boundaries for all the elements in $\Omega_N$ and $S(\Gamma):=\Pi_{T\in \Omega_N} L^2(\partial T)$ be the set of $L^2$ functions defined on $\Gamma$. For any $q \in S(\Gamma)$ and $\bq \in [S(\Gamma)]^d$,  we define their   averages $\{q\}, \{\bq\}$ and jumps $[q], [\bq]$ on the interior edges as follows. Suppose
$e$ is an interior edge shared by elements $T_+$ and $T_-$, we define the unit normal vectors $\bm{n}^+$ and  $\bm{n}^-$ on $e$ pointing exterior to $T_+$ and $T_-$, and then
\begin{flalign*}
[ q] \  =\  \, q^- \bm{n}^- \, +  q^+ \bm{n}^+, & \quad \{q\} = \frac{1}{2}( q^- + 	q^+), \\
[ \bq] \  =\  \, \bq^- \cdot \bm{n}^- \, +  \bq^+ \cdot \bm{n}^+, & \quad \{\bq\} = \frac{1}{2}( \bq^- + \bq^+).
\end{flalign*}

If $e$ is a boundary edge, then we let
$$[q] \  =q \bm{n}, \quad \{\bq\} =\bq,$$
where $\bm{n}$ is the outward unit normal.

Now we are ready to formulate the IPDG scheme for \eqref{elliptic eqn}. We look for $u_h \in \hat{\bV}_N^k $, such that
\begin{equation}
\label{sg}
B( u_h, v )= L(v),
 \quad  	\forall \,v \in \hat{\bV}_N^k
\end{equation}
where
\begin{flalign}
\label{eq:bilinear}
B(w, v )= &  \int_{\Omega} \bK \nabla w \cdot \nabla v\,d\bx
-  \sum_{\substack{e \in \Gamma}} \int_{e} \{\bK \nabla w\} \cdot [ v ]\,ds \notag \\
- & \sum_{\substack{e \in \Gamma}} \int_{e} \{\bK\nabla v\} \cdot [ w ]\,ds
+ \sum_{\substack{e \in \Gamma}}\frac{\sigma}{h} \int_{e} [w] \cdot [v]\,ds,
\end{flalign}
 and
\begin{flalign}
\label{eq:rhs}
L(v)= \int_\Omega f v \, d\bx-\int_{\partial \Omega} \left (\bK \nabla v \cdot \bm{n}+\frac{\sigma}{h} v \right ) g \, ds,
\end{flalign}
where $\sigma$ is a positive penalty parameter, $h=2^{-N}$ is the uniform mesh size in each dimension. The IP methods  developed in \cite{wheeler1978elliptic, arnold1982interior} are popular methods for elliptic equations, but usually with standard piecewise polynomial space. Here, by using the sparse finite element space, we can achieve a significant reduction in the size of the linear algebraic system \eqref{sg} especially when $d$ is large.

Several issues still have to be addressed to ensure real computational gains in implementation. First of all, it is obvious that we   need efficient algorithms to evaluate $L(v)$ when $v$ is taken to be the basis functions $v_{\bi, \bl}^{\bj}(\bx)$ in $\hat{\bV}_N^k.$ A standard integration scheme will incur computational complexity with exponential dependence on $d$. To avoid this, in our algorithm, we take full advantage of numerical integration on sparse grids developed in the literature \cite{smolyak1963quadrature, gerstner1998numerical}. In particular, for each basis $v_{\bi, \bl}^{\bj}(\bx)$,  its domain of dependence consists of $2^d$ smooth patches. Therefore, when evaluating
 $\int_\Omega f v_{\bi, \bl}^{\bj}\, d\bx$, we first divide this integral into $2^d$ parts accordingly, and then for each patch, we implement a sparse grid integration with Gauss quadrature, and the number of levels employed in this calculation is taken as $(I_{quad}-|\bl|_1/2)$, where $I_{quad}$ is a fixed integer chosen large enough to guarantee sufficient accuracy. The evaluation of $\int_{\partial \Omega} \left (\bK \nabla v \cdot \bm{n}+\frac{\sigma}{h} v \right ) g \, ds$ can be performed in a similar fashion.

  Another important factor when we assemble the linear system is the so-called \emph{unidirectional principle} \cite{bungartz1998finite}. To demonstrate the main ideas, we can see that evaluating $\int_\Omega \phi(\bx) \, d\bx$ with $\phi(\bx)=\phi_1(x_1) \ldots \phi_d(x_d)$ is equivalent to multiplication of one-dimensional integrals $\int_{[0,1]} \phi_1 \, dx_1 \cdots \int_{[0,1]} \phi_d \, dx_d.$ Therefore, computing  $\int_\Omega f v_{\bi, \bl}^{\bj}\, d\bx$ when  $f(\bx)$ is separable, i.e. $f(\bx)=f_1(x_1) \ldots f_d(x_d)$  or when $f(\bx)$ is a sum of separable functions is straightforward, because we only need some small overhead to compute one-dimensional integrals and assemble them to obtain the multi-dimensional integrations.

  The same discussion holds true for the computation of the bilinear term $B(u_h, v).$ For example, if we use a direct method, we need to evaluate $B(w, v)$ for $(w, v)$ being all possible basis functions in  $\hat{\bV}_N^k.$ From the definition of $B(w, v)$, each matrix element will involve four multi-dimensional integrations. If $\bK$ is separable (in particular, when $\bK$ is a constant function), due to the unidirectional principle, the matrices can be assembled fast.   When $\bK$ is a general function, we need to compute true high dimensional integrals. This difficulty is identified as one of the main challenging tasks for computing PDEs on sparse grids, see e.g. \cite{bungartz1998sparse, achatz2003higher}. In this case, we can either use the sparse grid integration procedure   mentioned above or by a computational procedure outlined as follows. Assume  $\bK$ to be a smooth function, then we can find $\bK_h = \hat{P}^{2k}_N \bK$ as the $L^2$ projection of $\bK$ onto the sparse finite element space $\hat{\bV}_N^{2k},$ and use $\bK_h$ in place of $\bK$ in the scheme. Notice that $\bK_h$ is a sum of separable functions, therefore the computation of the bilinear term is accelerated as the unidirectional case. The reason we use a higher order sparse finite element for projection of $\bK_h$ is to obtain exact evaluation of the volume integral $ \int_{\Omega} \bK \nabla u \cdot \nabla v\,d\bx$. However, this process does change the values of two 
  other terms up to approximation error of $\bK_h-\bK.$
Another aspect we did not explore is the efficient solver for the linear algebraic system. In the literature,  iterative methods
have been proposed based on the semi-coarsening approach and its sparse grid extensions \cite{mulder1989new, naik1993improved, griebel1990parallelizable, griebel1993multilevel, griebel1995tensor,pflaum1998multilevel, bungartz1997multigrid}. We leave those interesting implementation aspects to future study.



\subsection{Error Analysis}
\label{sec:error_analysis}

%
%
This section contains error estimates of the IPDG method on sparse grids.
Following \cite{arnold1982interior}, we define the  energy norm of a function $v \in H^2(\Omega_N)$ by
$$ ||| v ||| ^2 :=\sum_{\substack{T \in \Omega_N}}  \int_{T} |\nabla v|^2 \,d\bx \, + \sum_{\substack{e \in \Gamma}}h \int_{e} \left \{  \frac{\partial v}{\partial \bn} \right \}^2\,ds\, + \sum_{\substack{e \in \Gamma}}\frac{1}{h} \int_{e} [v]^2\,ds.$$

We review some basic properties of the bilinear operator $B(\cdot, \cdot).$
\begin{lem}[Orthogonality]
\label{lem:orth}
Let $u$ be the exact solution to \eqref{elliptic eqn}, and $u_h$ be the numerical solution to \eqref{sg}, then
$$B( u-u_h, v )=0,   \quad \forall  v \in \hat{\bV}_N^k.$$
\end{lem}
\noindent{\it Proof}:
Using integration by parts, we can easily show $B(u, v)=0,  \, \forall  v \in \hat{\bV}_N^k.$ The Galerkin orthogonality immediately follows.
 \hfill $ \blacksquare $

\begin{lem}[Boundedness\cite{arnold1982interior, Arnold_2002_SIAM_DG}]
\label{lem:bound}
There exists a positive constant $C_b$, depending only on $K_1, \sigma$, such that
$$B( w, v ) \le C_b |||w |||\cdot ||| v |||,   \quad \forall \,w, v \in H^2(\Omega_N).$$
\end{lem}

\begin{lem}[Stability\cite{arnold1982interior,Arnold_2002_SIAM_DG}]
\label{lem:stab}
When $\sigma$ is taken large enough, there exists a positive constant $C_s$ depending only on $K_0$, such that
$$B( v, v ) \ge C_s  ||| v |||^2,   \quad \forall \,  v \in \hat{\bV}_N^k.$$
\end{lem}

\begin{thm}
Let $u$ be the exact solution to \eqref{elliptic eqn}, and $u_h$ be the numerical solution to \eqref{sg}. For  $k \geq 1$,  $u \in \mathcal{H}^{p+1}(\Omega)$,  $1 \leq t \leq \min \{p, k\}$,
 $N\geq 1$, $d \geq 2$,  we have
$$
|||u-u_h|||    \leq \left (1+\frac{C_b}{C_s} \right) \sqrt{C_d} \underline{c}_{k,t}  C^* 2^{-Nt} |u|_{\mathcal{H}^{t+1}(\Omega)}.
$$
where $C_d$ is a constant that depends on $d$ linearly. $C_b, C_s$ are defined in Lemmas \ref{lem:bound} and \ref{lem:stab}.
$$
C^*=\max \left(\sqrt{d^2 \kappa_0^{2d-2}+3 d^3  \kappa_1^{2d}  + 2d^4 \kappa_0^{2d-2} 2^{-2N} }, \sqrt{d^2 \kappa_0^{2d-2}+d^3  \kappa_1^{2d} +2d+ 2d^4 \kappa_0^{2d-2} 2^{-2N}} \right),
$$
where $\kappa_s=\kappa_s(k, t,  N), \,s=0, 1$ and $\underline{c}_{k,t}$ are defined in Theorem \ref{thm:appx}.
\end{thm}
\noindent{\it Proof}:
Choose any function $u_I \in \hat{\bV}_N^k$, then we decompose the error into $e=u-u_h=(u - u_I)+(u_I - u_h)$.
By C\'{e}a's lemma, using Lemma \ref{lem:orth}, \ref{lem:bound} and \ref{lem:stab},
$$
C_s |||u_I-u_h|||^2 \le  B_h(u_I-u_h, u_I-u_h)  = B_h(u_I-u, u_I-u_h) \leq C_b  |||u-u_I |||\cdot |||u_I-u_h|||.
$$
Therefore, $|||u_I-u_h|| \leq \frac{C_b}{C_s} |||u-u_I |||,$ and
$$
|||e||| \leq |||u_I-u_h|||+|||u-u_I||| \leq \left (1+\frac{C_b}{C_s} \right) |||u-u_I |||.
$$
We have
$$
|||e||| \leq  \left (1+\frac{C_b}{C_s} \right)  \inf_{u_I \in  \hat{\bV}_N^k}|||u-u_I ||| \leq \left (1+\frac{C_b}{C_s} \right)  |||u- \underline{\hat{\bP}}_N^k u|||,
$$
where the projection operator $\underline{\hat{\bP}}_N^k$ has been specified in
Theorem \ref{thm:appx}.

Next, we need to bound  the energy norm of $\eta:=u- \underline{\hat{\bP}}_N^k u$. Recall that
$$ ||| \eta ||| ^2 =\sum_{\substack{T \in \Omega_N}}  \int_{T} |\nabla \eta|^2 \,d\bx \, + \sum_{\substack{e \in \Gamma}}h \int_{e} \left \{  \frac{\partial \eta}{\partial \bn} \right \}^2\,ds\, + \sum_{\substack{e \in \Gamma}}\frac{1}{h} \int_{e} [\eta]^2\,ds.$$

The first term in the summation is
$$\sum_{\substack{T \in \Omega_N}}  \int_{T} |\nabla \eta|^2 \,d\bx = |\eta|^2_{H^1(\Omega_N)}.$$

To bound the second and third terms, we use the trace inequalities \cite{arnold1982interior}:
\begin{flalign*}
\left | \left |\phi \right |\right|^2_{L^2(\partial T)} \le C_d \left(\frac{1}{h} ||\phi||^2_{L^2(T)} + h |\phi|^2_{H^1(T)}  \right), \quad & \forall \, \phi \in H^1(T) \\
\left | \left |\frac{\partial \phi}{\partial \bn} \right |\right|^2_{L^2(\partial T)} \le C_d \left(\frac{1}{h} |\phi|^2_{H^1(T)} + h |\phi|^2_{H^2(T)}  \right), \quad & \forall \, \phi \in H^2(T)
\end{flalign*}
where $C_d$ is a generic constant that depends on $d$ linearly. It may have different values in each occurrence in the proof. Sum over all the elements, we get
$$
\sum_{\substack{e \in \Gamma}}h \int_{e} \left \{  \frac{\partial \eta}{\partial \bn} \right \}^2\,ds \le C_d \left ( |\eta|^2_{H^1(\Omega_N)}+ h^2 |\eta|^2_{H^2(\Omega_N)} \right)
$$
$$
 \sum_{\substack{e \in \Gamma}}\frac{1}{h} \int_{e} [\eta]^2\,ds \le C_d \left ( \frac{1}{h^2} ||\eta||^2_{L^2(\Omega_N)}+  |\eta|^2_{H^1(\Omega_N)} \right).
$$
In summary,
$$
||| \eta ||| ^2  \le C_d \left ( \frac{1}{h^2} ||\eta||^2_{L^2(\Omega_N)}+  |\eta|^2_{H^1(\Omega_N)} +h^2 |\eta|^2_{H^2(\Omega_N)} \right).
$$
By Theorem \ref{thm:appx},
\begin{flalign*}
||| \eta ||| ^2  \le C_d \, \underline{c}_{k,t}^2\left ( d^2 \kappa_0^{2d-2}+3 d^3  \kappa_1^{2d}  + 2d^4 \kappa_0^{2d-2} 2^{-2N}        \right ) 2^{-2N t} |u|^2_{\mathcal{H}^{t+1}(\Omega)}, \quad k \geq 2 \\
||| \eta ||| ^2  \le C_d \, \underline{c}_{k,t}^2\left ( d^2 \kappa_0^{2d-2}+d^3  \kappa_1^{2d} +2d+ 2d^4 \kappa_0^{2d-2} 2^{-2N}        \right ) 2^{-2N t} |u|^2_{\mathcal{H}^{t+1}(\Omega)}, \quad k =1
\end{flalign*}
where we have used the shorthand notation $\kappa_s=\kappa_s(k, t,  N), \,s=0, 1.$ Let's define
$$
C^*=\max \left(\sqrt{d^2 \kappa_0^{2d-2}+3 d^3  \kappa_1^{2d}  + 2d^4 \kappa_0^{2d-2} 2^{-2N} }, \sqrt{d^2 \kappa_0^{2d-2}+d^3  \kappa_1^{2d} +2d+ 2d^4 \kappa_0^{2d-2} 2^{-2N}} \right).
$$
Then,
$$
||| \eta ||| \le \sqrt{C_d} \underline{c}_{k,t}  C^* 2^{-Nt} |u|_{\mathcal{H}^{t+1}(\Omega)}.
$$
Therefore, we have proved the error estimate in the energy norm
$$
|||e|||    \leq \left (1+\frac{C_b}{C_s} \right) \sqrt{C_d} \underline{c}_{k,t}  C^* 2^{-Nt} |u|_{\mathcal{H}^{t+1}(\Omega)}.
$$
 \hfill $ \blacksquare $

This theorem implies a convergence rate of $O(h^{k})$ up to the polylogarithmic term $|\log_2 h|^{d-1}$ in the energy norm when $u$ is smooth enough. However, we do require more regularity of $u$ compared with IPDG methods using traditional piecewise polynomial space.

 \begin{rem}
Proving convergence in $L^2$ norm with the standard duality argument will encounter some difficulties in this framework. For example, let $\varphi$ be the solution to the adjoint problem
$$
-\nabla \cdot (\bK \nabla \varphi)=u-u_h \quad \textrm{on} \,\, \Omega, \qquad \varphi=0 \quad \textrm{on} \,\, \partial \Omega.
$$
Since $\Omega$ is convex, we have $\varphi \in H^2(\Omega)$, and $||\varphi||_{H^2(\Omega)} \leq C ||u-u_h||_{L^2(\Omega)}.$ From the adjoint consistency of the IP method, we get
$$
||u-u_h||^2_{L^2(\Omega)}=B(u-u_h, \varphi).
$$
However, to proceed from here, we will need to define an interpolant of $ \varphi$: $\varphi_I \in \hat{\bV}_N^1$, and bound $|||\varphi-\varphi_I|||$. From our previous discussion, this will require a bound in $||\varphi||_{\mathcal{H}^{2}(\Omega)}.$ This is a stronger norm than the classical $H^2$ norm, and cannot be  bounded by $||u-u_h||^2_{L^2(\Omega)}$.
 \end{rem}

\section{Numerical Results}
\label{sec: SGDG numerical}
In this section, we provide multi-dimensional numerical results to demonstrate the performance of our sparse grid IPDG scheme.

\subsection{Two-dimensional Results}

In this subsection, we gather computational results for two-dimensional case. The penalty constant in this subsection is chosen to be $\sigma = 10$ for $k=1,$ and $\sigma = 20$ for $k=2$.

\begin{figure}[htp]
  \centering
  \begin{tabular}{c c}

(a)\includegraphics[width=.48\textwidth]{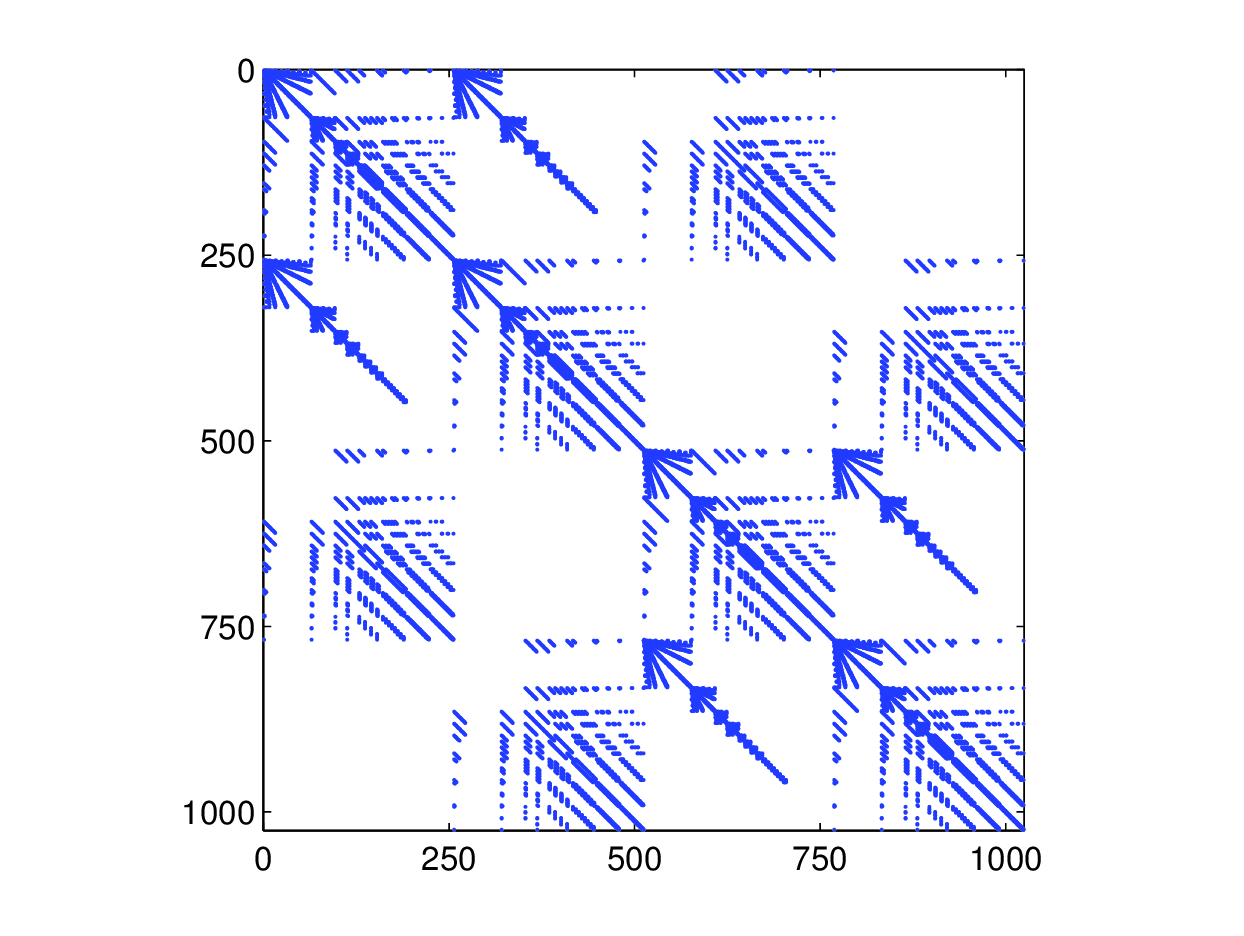}
(b)\includegraphics[width=.48\textwidth]{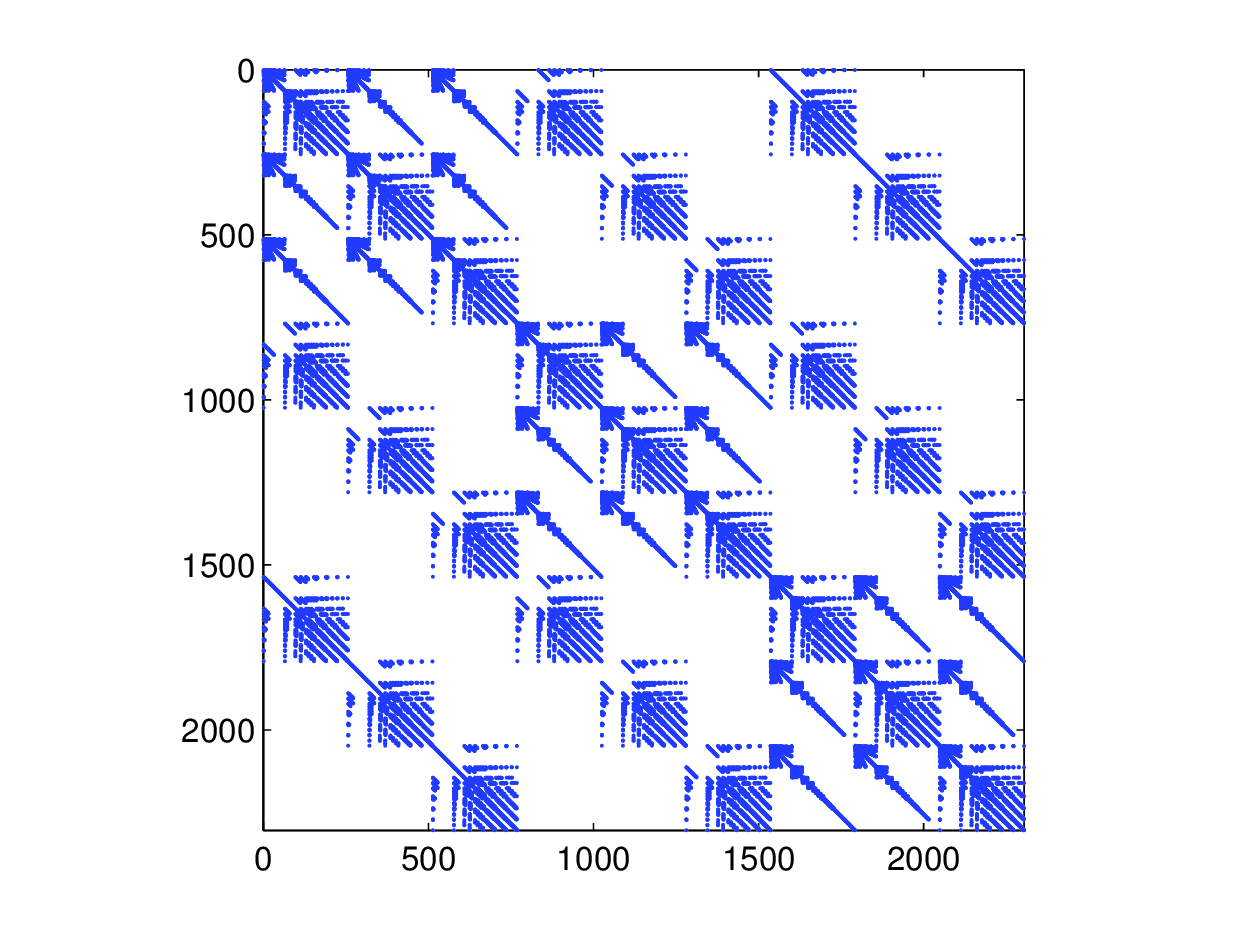}
\end{tabular}
  \caption{Example \ref{ex:const_2D}.  The sparsity patterns of the matrices computed by the space $\hat{\bV}_6^k$ with $ k=1, 2$ and $N = 6$. Each dot represents a non-zero element in the stiffness matrix. (a): $k=1$, (b): $k=2.$}
 \label{fig:2Dconst_spy}
\end{figure}

\begin{exa}\rm We solve the following two-dimensional problem with constant   coefficient on $\Omega = [0,1]^2 $. \label{ex:const_2D}
\begin{align}
\begin{cases}
- \Delta u = 0, &  \qquad \bx \in \Omega ,\\
u = \sin(\pi x_1)\, \frac{\sinh(\pi x_2)}{\sinh(\pi)}, & \qquad \bx \in \partial \Omega ,
  \end{cases}
\end{align}
for which the exact solution is $u = \sin(\pi x_1)\, \frac{\sinh(\pi x_2)}{\sinh(\pi)}$. We test the scheme with  $k=1$ and $k=2$ on different levels of meshes.
 Numerical errors and orders of accuracy measured in $L^1, L^2, L^\infty$ and $H^1$ norms are listed in Table \ref{table:2DconstBC}. We observe $k$-th order of accuracy for $H^1$ norm, close to $(k+1)$-th order accuracy for $L^1, L^2$ norms and slightly less than $(k+1)$-th order accuracy for the $L^\infty$ norm. The results agree with the error estimates performed in the previous section. The slight order reduction for $L^1, L^2$, and $H^1$ norms is similar to the one observed in Section \ref{sec:approx}.

In addition, we provide the sparsity patterns and the condition number of the stiffness matrices for $k=1$ and $k=2$   in Figure \ref{fig:2Dconst_spy} and Table \ref{table:2Dconst_spy}. From Table \ref{table:2Dconst_spy}, the number of nonzero elements scales as $O(SGDOF^{1.5})$, where $SGDOF$ is the number of degrees of freedom of the space used. This is a denser matrix than the one generated from traditional piecewise polynomial space, for which one element can only interact with itself and its immediate neighbors. For sparse grids, the basis functions in $\hat{\bV}_N^k$ are no longer locally defined due to the hierarchical structure. However, fixing any point $\bx$ in the domain $\Omega$, the number of bases with nonzero values at $\bx$ is a constant for each level. In this sense, the interaction between the bases is still limited when compared with a global approximation scheme.

\begin{table}
\caption{Numerical errors and orders of accuracy for Example \ref{ex:const_2D} computed by the  space $\hat{\bV}_N^k$ with $k=1, 2$ and indicated $N$. 
}
\vspace{2 mm}
\centering
\begin{tabular}{ c c c c c c c  c c }
\hline
N &$L^1$ error & order & $L^2$ error & order & $L^\infty$ error & order  &$H^1$ error & order \\
\hline
 &     & &  & $ k=1$    & &  &  & \\
\hline
 3          &4.49E-03  &          &  6.97E-03 &          &  3.26E-02 &            &1.77E-01 & \\
4             &1.18E-03  &  1.93 & 1.93E-03  &  1.85 & 9.71E-03  &1.75 & 8.80E-02& 1.01 \\
5           &   3.03E-04 &1.96  &5.09E-04&  1.92 &  3.19E-03 &  1.60   &4.36E-02 & 1.01 \\
6           & 7.68E-05  & 1.98 & 1.32E-04 &  1.98 &  9.68E-04 &1.72  &2.16E-02 & 1.01 \\
\hline
 &    & & &  $k=2$   & &  &  &   \\
 \hline
3              &9.52E-05 &         &  1.33E-04&          & 5.74E-04 &          &7.61E-03&     \\
4               &1.42E-05   & 2.75 &    2.03E-05&  2.71 & 9.65E-05 & 2.57&  1.91E-03 & 1.99 \\
5             &2.05E-06  &2.79 & 3.02E-06 &  2.75 & 1.59E-05 &  2.60 &  4.78E-04 & 2.00\\
6          &   2.89E-07  & 2.83 & 4.36E-07&  2.79 & 2.66E-06 &  2.58 &1.19E-04 & 2.00 \\
\hline

\end{tabular}
\label{table:2DconstBC}
\end{table}
\end{exa}

\begin{table}
\caption{Sparsity and condition number of the stiffness matrix. Example \ref{ex:const_2D} computed by the space $\hat{\bV}_N^k$ with $ k=1, 2$. SGDOF is the number of degrees of freedom used for the sparse grid DG scheme. 
NNZ is the number of nonzero elements  in the stiffness matrix. Order=$\log(\textrm{NNZ})/\log(\textrm{SGDOF}).$}
\vspace{2 mm}
\centering
\begin{tabular}{c c c c  c }
\hline
N & SGDOF & NNZ &   Order & Condition Number \\
\hline
&  & & $k=1$      &\\
\hline
 3  &   80    & 992 &    1.57  &  3.58E+02 \\
4  &   192   &  3216 &  1.54  & 1.43E+03    \\
5   &   448    & 9168 &        1.49  & 5.68E+03\\
6   &   1024   & 24144 &      1.45 & 2.26E+04\\
\hline
 &  & & $k=2$ & \\
 \hline
3  & 180    & 3456 &   1.57  &  1.40E+03   \\
4    & 432     &11124 &    1.54 & 5.49E+03\\
5  & 1008    &31596 &  1.50 & 2.16E+04\\
6  &2304     &83028 &    1.46 & 8.58E+04\\
\hline

\end{tabular}
\label{table:2Dconst_spy}
\end{table}

\begin{exa}\rm We solve the following two-dimensional problem with smooth variable coefficient on $\Omega = [0,1]^2$.
\label{ex:smth}
\begin{align}
\begin{cases}
- \nabla \cdot ((\sin(x_1 \,x_2)  +1 ) \,  \nabla u) = f , &  \qquad \bx \in \Omega, \\
u = 0, & \qquad \bx \in \partial \Omega.
 \end{cases}&
\end{align}
$f$ is a given function such that the exact solution is $u = \sin(\pi x_1)\, \sin(\pi x_2)$. As discussed in previous Section, when assembling the stiffness matrix, we first obtain the $L^2$ projection of coefficient $K$ in space $\hat{\mathbf{V}}_N^{2k}$, then make use of the unidirectional principle to save computational cost. The numerical results are provided in Table \ref{table:smth}. The conclusion is very similar to the constant coefficient case Example \ref{ex:const_2D}. We observe $k$-th order of accuracy for $H^1$ norm, close to $(k+1)$-th order of  accuracy for $L^1, L^2$ norms and slightly less than $(k+1)$-th order of accuracy for the $L^\infty$ norm.
\begin{table}
\caption{Numerical errors and orders of accuracy for Example \ref{ex:smth} computed by the space $\hat{\bV}_N^k$ with $ k=1, 2$ and indicated $N$. 
}
\vspace{2 mm}
\centering
\begin{tabular}{c  c c c c c c c c }
\hline
N  &$L^1$ error & order & $L^2$ error & order & $L^\infty$ error & order  &$H^1$ error & order \\
\hline
  &   &  &   & $ k=1$    & &   &  & \\
\hline
 3              &1.30E-02  &          &  1.65E-03 &          &  4.50E-02 &            &3.37E-01 & \\
4               &3.18E-03  &  2.03 & 4.08E-03  &  2.01 & 1.34E-02  &1.74 & 1.66E-01& 1.02 \\
5              &   7.81E-04 & 2.02  &1.01E-03&  2.01 &  4.43E-03 &  1.60   &8.26E-02 & 1.01 \\
6             & 1.94E-04  & 2.01 & 2.55E-04 &  1.99 &  1.48E-03&  1.58   &4.11E-02 & 1.00 \\
\hline
 &    & & &   $k=2$    &  &  & &   \\
 \hline
3          &1.77E-04 &         &  2.17E-04&          &  5.74E-04 &          &1.35E-02&     \\
4           &2.71E-05   & 2.70 &    3.37E-05&  2.69 &1.01E-04& 2.51&  3.37E-03 & 2.00 \\
5          & 3.99E-06  &2.76 & 5.08E-06 &  2.73 & 1.91E-05&  2.40 &   8.41E-04 & 2.00\\
6         &   5.67E-07  & 2.82 &7.37E-07&  2.78 & 2.99E-06 &  2.67 &2.10E-04 & 2.00 \\
\hline

\end{tabular}
\label{table:smth}

\end{table}

\end{exa}

\begin{exa}\rm We solve the following two-dimensional problem with discontinuous coefficient on $\Omega = [0,1]^2 $.
\label{ex:nonsmth_2D}
\begin{align}
\ \ \ \ \ &\begin{cases}
- \nabla \cdot ((    \sign( (x_1-0.5)(x_2-0.5) ) +2   )\,  \nabla u) = f,  & \qquad \bx \in \Omega,\\
 u = 0, & \qquad \bx \in \partial \Omega.
  \end{cases}&
\end{align}
$f$ is a given function such that the exact solution is $u = \sin(\pi x_1)\, \sin(\pi x_2)$.    In this example, the coefficient has four constant panels on $\Omega$, but the solution itself is smooth.  We want to use this example to demonstrate the performance of the proposed scheme when solving problems with discontinuous coefficient. We provide the numerical results with $k=1$ and $k=2$ in Table \ref{table:nonsmth_2D}.   Notice that in this example, the discontinuity of the coefficient is along the cell boundaries for all our computational grids, and the convergence rates are similar to the previous two examples.

\begin{table}
\caption{Numerical errors and orders of accuracy for Example \ref{ex:nonsmth_2D} computed by the space $\hat{\bV}_N^k$ with $ k=1, 2$   and indicated $N$.
}
\vspace{2 mm}
\centering
\begin{tabular}{c c c c c c c c c }
\hline
N &$L^1$ error & order & $L^2$ error & order & $L^\infty$ error & order  &$H^1$ error & order \\
\hline
  &   &   &  & $ k=1$    & &  &  & \\
\hline
 3       &1.24E-02  &          &  1.57E-02 &          &  4.55E-02 &            &3.33E-01 & \\
4        &3.07E-03  &  2.02 &3.94E-03  &  2.00 & 1.36E-02  &1.75 & 1.66E-01& 1.00 \\
5        &  7.58E-04 &2.02 &9.78E-04& 2.01 &  4.50E-03 &  1.59  &8.32E-02 & 1.00 \\
6       & 1.89E-04  & 2.00 & 2.46E-04 &  1.99 &  1.50E-03 &1.58  &4.16E-02 & 1.00 \\
\hline
 &    & &  & $k=2$   &  &  & &   \\
 \hline
3        &1.96E-04 &         &  2.59E-04&          & 1.21E-03 &          &1.56E-02&     \\
4        &2.72E-05   & 2.85 &    3.50E-05&  2.89 &1.55E-04& 2.96&  3.70E-03 & 2.08 \\
5        &3.85E-06  &2.82 & 4.94E-06 &  2.82 & 2.05E-05 &  2.92 &  8.93E-04 & 2.05\\
6        &   5.36E-07  & 2.84 &7.02E-07&  2.82 & 3.01E-06 &  2.82 &2.19E-04 & 2.03 \\
\hline

\end{tabular}
\label{table:nonsmth_2D}
\end{table}

\end{exa}

\subsection{Three-dimensional Results}

In this subsection, we gather computational results for three-dimensional elliptic equations. The penalty constant in this subsection is chosen to be $\sigma = 15$ for  $k=1,$ and $\sigma = 30$ for $k=2$.

\begin{figure}[htp]
  \centering
  \begin{tabular}{c c}

(a)\includegraphics[width=.47\textwidth]{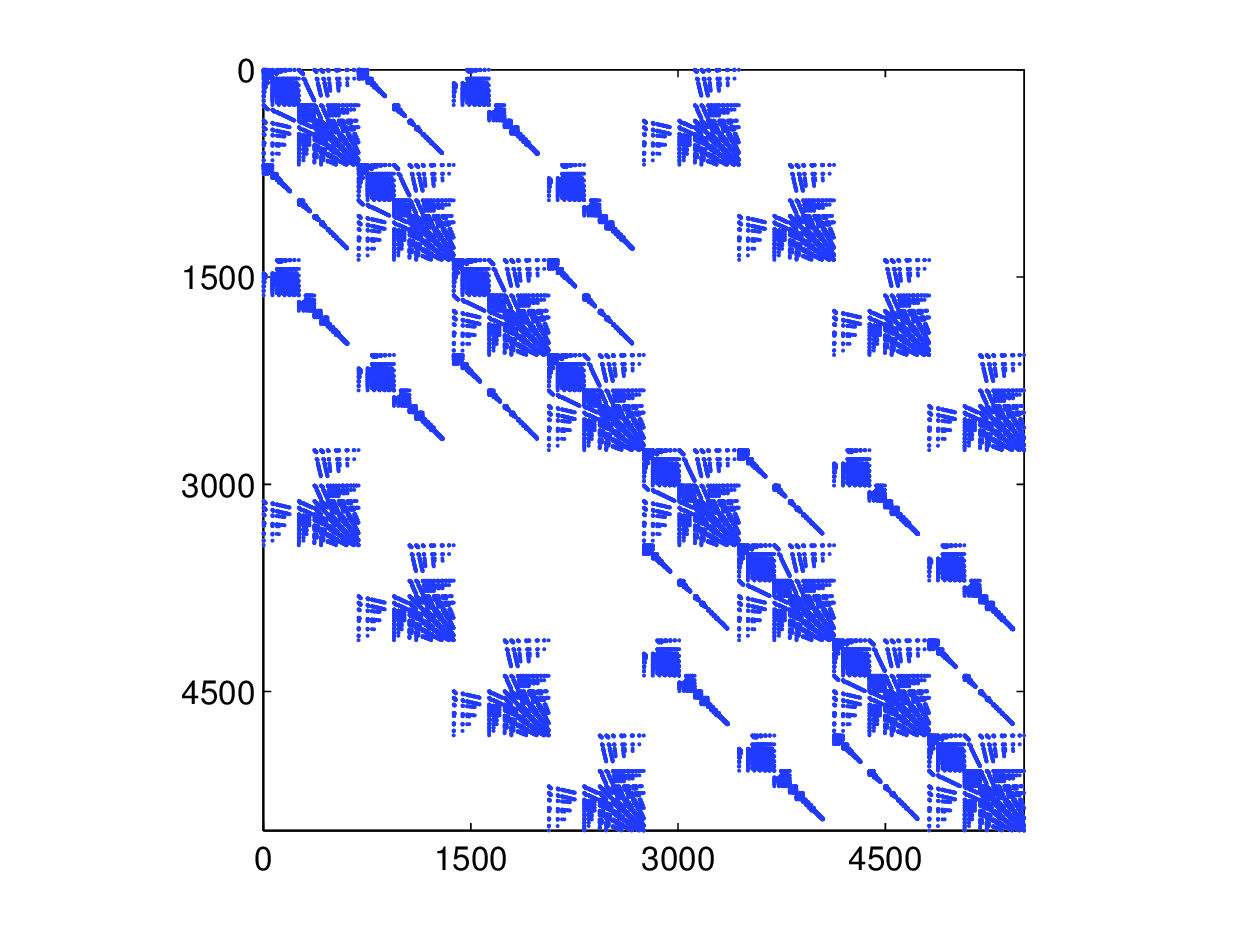}
(b)\includegraphics[width=.47\textwidth]{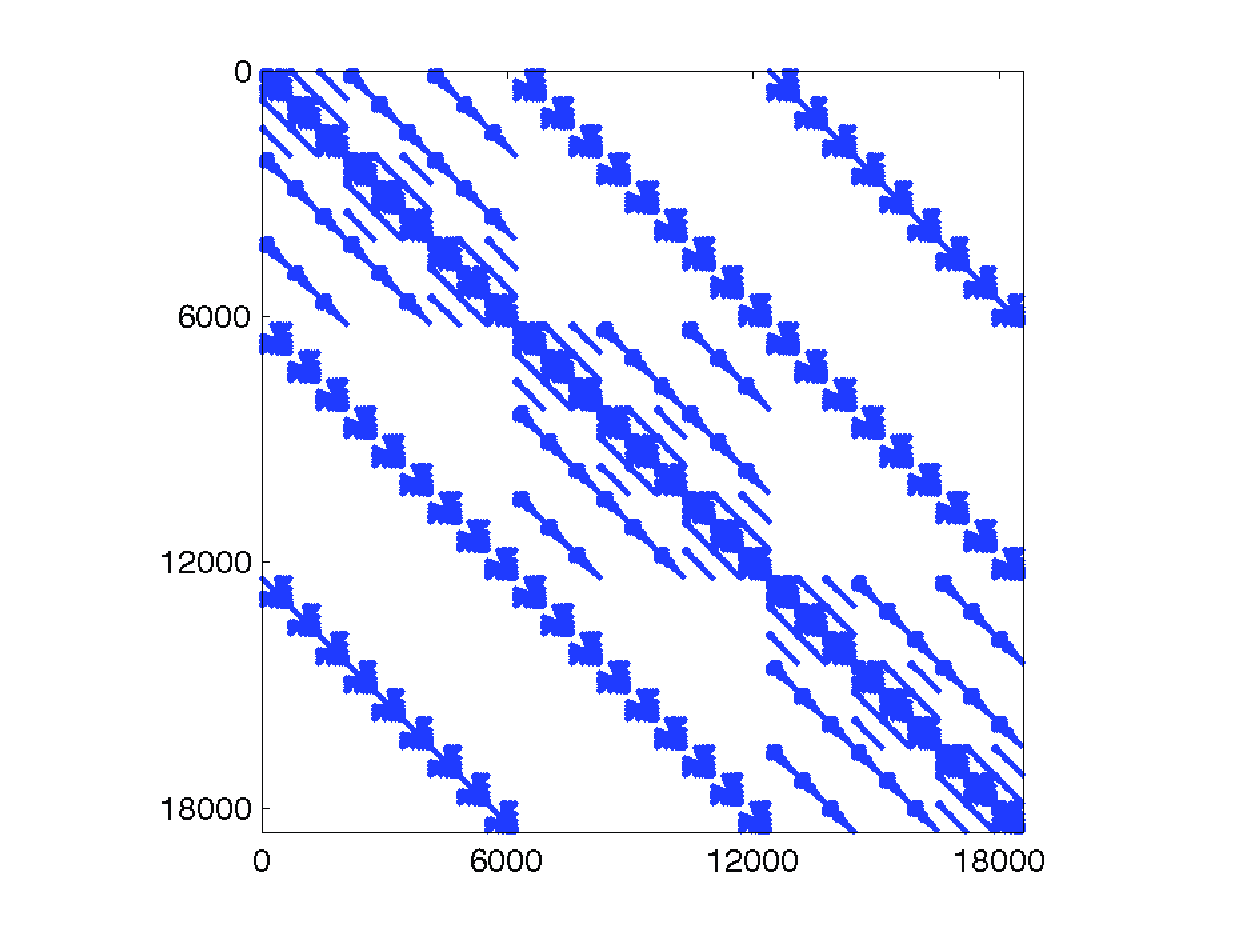}
\end{tabular}
  \caption{Example \ref{ex:const_3D}.  The sparsity patterns of the matrices computed by the space $\hat{\bV}_6^k$ with (a) $k=1$ and (b) $k=2$.   Each dot represents a non-zero element in the stiffness matrix.}
 \label{fig:3Dconst_spy}
\end{figure}

\begin{exa}\rm We solve the following three-dimensional problem with constant  coefficient on $\Omega = [0,1]^3 $. \label{ex:const_3D}
\begin{align}
\begin{cases}
- \Delta u = 0,   &\qquad \bx \in \Omega,\\
u = \sin(\pi x_1)\, \sin(\pi x_2)\, \frac{\sinh(\sqrt{2} \pi x_3)}{\sinh(\sqrt{2} \pi)},  &\qquad \bx \in \partial \Omega,
  \end{cases}&
\end{align}
where the exact solution is $u = \sin(\pi x_1)\, \sin(\pi x_2)\, \frac{\sinh(\sqrt{2} \pi x_3)}{\sinh(\sqrt{2} \pi)}$.
We provide the numerical results with $k=1$ and $k=2$ in Table \ref{table:3DconstBC}.  For this three-dimensional example, we obtain $k$-th order for the $H^1$ norm, and close to $(k+1)$-th order for the $L^1$ and $L^2$ norm, but $L^\infty$ error seems to degrade to $(k+\frac{1}{2})$-th order. However, upon closer examination, similar behaviors have been observed in Section \ref{sec:approx} for the $L^2$ projection error of a smooth function onto $\hat{\bV}_N^k.$
The sparsity patterns of the stiffness matrices for $k=1$ and $k=2$ are reported in Figure \ref{fig:3Dconst_spy} and Table \ref{table:3Dconst_spy}. From Table \ref{table:3Dconst_spy}, we can see the stiffness matrices scale less than the two-dimensional examples, i.e., the number of nonzero elements scales as $O(SGDOF^{1.4})$, where $SGDOF$ is the number of degrees of freedom of the space used.

\begin{table}
\caption{Numerical errors and orders of accuracy for Example \ref{ex:const_3D} computed by the space $\hat{\bV}_N^k$ with $ k=1, 2$ and indicated $N$. 
}
\vspace{2 mm}
\centering
\begin{tabular}{c c c c c c c  c c }
\hline
N & $L^1$ error & order & $L^2$ error & order & $L^\infty$ error & order  &$H^1$ error & order \\
\hline

  &     & &  & $ k=1$    & &  &  & \\
\hline
 3   &1.29E-02  &          &  2.19E-02 &          &  1.09E-01 &            &2.85E-01 & \\
4    &4.05E-03  &  1.67 & 6.98E-03  &  1.65 & 4.75E-02  &1.20 & 1.44E-01& 0.98 \\
5   &   1.07E-03 &1.92  &1.94E-03&  1.85 &  2.34E-02 &  1.02   &7.02E-02 & 1.04 \\
6   & 2.76E-04  & 1.96 &5.22E-04 &  1.89 &  8.44E-03 &1.47  &3.39E-02 & 1.05 \\
\hline
 &    & &  & $k=2$   &  &  & &   \\
 \hline
3   &1.41E-04 &         &  2.06E-04&          & 1.26E-03 &          &1.05E-02&     \\
4   &2.51E-05   & 2.49 &    3.80E-05&  2.44 &3.35E-04& 1.91&  2.72E-03 & 1.95 \\
5   &4.18E-06  &2.59 & 6.49E-06 & 2.55 & 6.51E-05 &  2.36 &  6.87E-04 & 1.98\\
6  &  6.69E-07  & 2.64 &1.06E-06&  2.62 &1.09E-05&  2.58 &1.72E-04 & 2.00 \\
\hline

\end{tabular}
\label{table:3DconstBC}
\end{table}
\end{exa}

\begin{table}
\caption{Sparsity and condition number of the stiffness matrix. Example \ref{ex:const_3D} computed by the space $\hat{\bV}_N^k$ with $ k=1, 2$. SGDOF is the number of degrees of freedom used for the sparse grid DG scheme. 
NNZ is the number of nonzero elements  in the stiffness matrix. Order=$\log(\textrm{NNZ})/\log(\textrm{SGDOF}).$}
\vspace{2 mm}
\centering
\begin{tabular}{c c c c c }
\hline
N &SGDOF & NNZ  & Order & Condition Number \\
\hline

&  & & $ k=1$   &\\
\hline
 3 &   304  &3760  &   1.43  &3.73E+02 \\
4  &   832   & 14080  & 1.42  &1.51E+03   \\
5 &  2176 & 45760 &    1.39 & 5.97E+03  \\
6 &   5504   & 135872 &   1.37 & 2.36E+04\\
\hline
 &  & & $k=2$     & \\
 \hline
3    & 1026   &20250      & 1.43  &1.58E+03    \\
4    &2808  &74628 & 1.41 & 5.98E+03\\
5   & 7344   &240516 & 1.39 &2.32E+04 \\
6  &18576     &710532 & 1.37 & 9.15E+04 \\
\hline

\end{tabular}
\label{table:3Dconst_spy}
\end{table}

\begin{exa}\rm We solve the following three-dimensional problem with smooth variable coefficient on $\Omega = [0,1]^3$.
\label{ex:3Dsmth}
\begin{align}
\begin{cases}
- \nabla \cdot ((\sin(x_1 x_2 x_3) +1 )\,  \nabla u) = f , &  \qquad \bx \in \Omega,\\
u = 0, &  \qquad \bx \in \partial\Omega,
  \end{cases}&
\end{align}
$f$ is a given function such that the exact solution is $u = \sin(\pi x_1)\, \sin(\pi x_2)\, \sin(\pi x_3)$. Again, for such a non-separable coefficient $K$, the $L^2$ projection onto space $\hat{\mathbf{V}}_N^{2k}$ is obtained and the unidirectional principle is applied when we assemble the stiffness matrix.  
We provide the numerical results for $k=1$ and $k=2$  in Table \ref{table:3Dsmth}. The conclusion of Example \ref{ex:const_3D} holds here, including the comparable convergence rate and similar order reduction of the $L^1$, $L^2$, and $L^\infty$ errors.

\begin{table}
\caption{Numerical errors and   orders of accuracy for Example \ref{ex:3Dsmth}  computed by the space $\hat{\bV}_N^k$ with $ k=1, 2$ and indicated $N$. 
}
\centering
\begin{tabular}{c c c c c c c c c }
\hline
N &$L^1$ error & order & $L^2$ error & order & $L^\infty$ error & order  &$H^1$ error & order \\
\hline

  &   &  & &   $ k=1$    & &  &  & \\
\hline
 3           &2.64E-02  &          &  3.40E-02 &          &  1.55E-01 &            &4.32E-01 & \\
4               &6.23E-03  &  2.08 &8.58E-03  &  1.98 & 3.54E-02  &2.13 & 2.04E-01& 1.08 \\
5            &  1.49E-03 &2.06 & 2.10E-03 & 2.03 &  2.07E-02 &  0.77  &9.82E-02 & 1.06 \\
6          & 3.68E-04  & 2.02 & 5.32E-04 &  1.98 &  7.58E-03 &1.45  &4.80E-02 & 1.03 \\
\hline
 &    & &   & $k=2$   & &  & &   \\
 \hline
 3       &1.63E-04&         & 2.05E-04&          & 8.24E-04&          &1.19E-02&     \\
4       &2.88E-05  & 2.50 &    3.66E-05& 2.48 &1.63E-04&2.34&  3.00E-03 & 1.98 \\
5    &4.72E-06  &2.61 & 6.06E-06 &  2.60 &2.73E-05 &  2.58 & 7.54E-04 & 2.00\\
6    &7.42E-07 &2.67 &9.58E-07 &  2.66 &5.80E-06 &  2.23 &1.88E-04 & 2.00\\
\hline

\end{tabular}
\label{table:3Dsmth}
\end{table}

\end{exa}

\newpage
 \subsection{Four-dimensional Results}
 In this subsection, we gather computational results for four-dimensional case. The penalty constant in this subsection is chosen to be $\sigma = 30$ for $k=1,$ and $\sigma = 60$ for $k=2$.

\begin{exa}\rm We solve the following four-dimensional problem with constant variable coefficient on $\Omega = [0,1]^4 $.
\label{ex:const_4D}
\begin{align}
\begin{cases}
- \Delta u = 0  &  \qquad \bx \in \Omega,\\
u = \sin(\pi x_1)\, \sin(\pi x_2)\, \sin(\pi x_3)  \, \frac{\sinh(\sqrt{3} \pi x_4)}{\sinh(\sqrt{3} \pi)}  &\qquad \bx \in \partial \Omega
 \end{cases}&
\end{align}
where the exact solution is $u = \sin(\pi x_1)\, \sin(\pi x_2)\, \sin(\pi x_3)  \, \frac{\sinh(\sqrt{3} \pi x_4)}{\sinh(\sqrt{3} \pi)} $.

We list numerical errors and orders in Table \ref{table:4DconstBC}. 
For this four-dimensional example, we obtain $k$-th order for the $H^1$ norm, and close to $(k+\frac{1}{2})$-th order for the $L^1$ and $L^2$ norms. The $L^\infty$ error seems to fluctuate. However,  we expect the order of accuracy in $L^\infty$ norm to grow upon refinement.  The numerical results when fixing two of the coordinates are plotted in Figure \ref{fig:const_4D} for $k=1$ and $k=2$ when $N=4$.  The left  and right figures are obtained with degrees of freedom equal to 3072 and 15552, respectively. We can see the higher order scheme offers a  much more accurate solution than the lower order scheme, and this is also verified by Table \ref{table:4DconstBC}.
The sparsity patterns and condition number of the stiffness matrices for $k=1$ and $k=2$ are reported in Figure \ref{fig:4Dconst_spy} and Table \ref{table:4Dconst_spy}. From Table \ref{table:4Dconst_spy}, we can see the stiffness matrices scale less than the two-dimensional and three-dimensional examples, i.e., the number of nonzero elements scales like $O(SGDOF^{1.35})$, where $SGDOF$ is the degree of freedom of the space used.


\begin{figure}[htp]
  \centering
  \begin{tabular}{c c}

(a) \includegraphics[width=.5\textwidth]{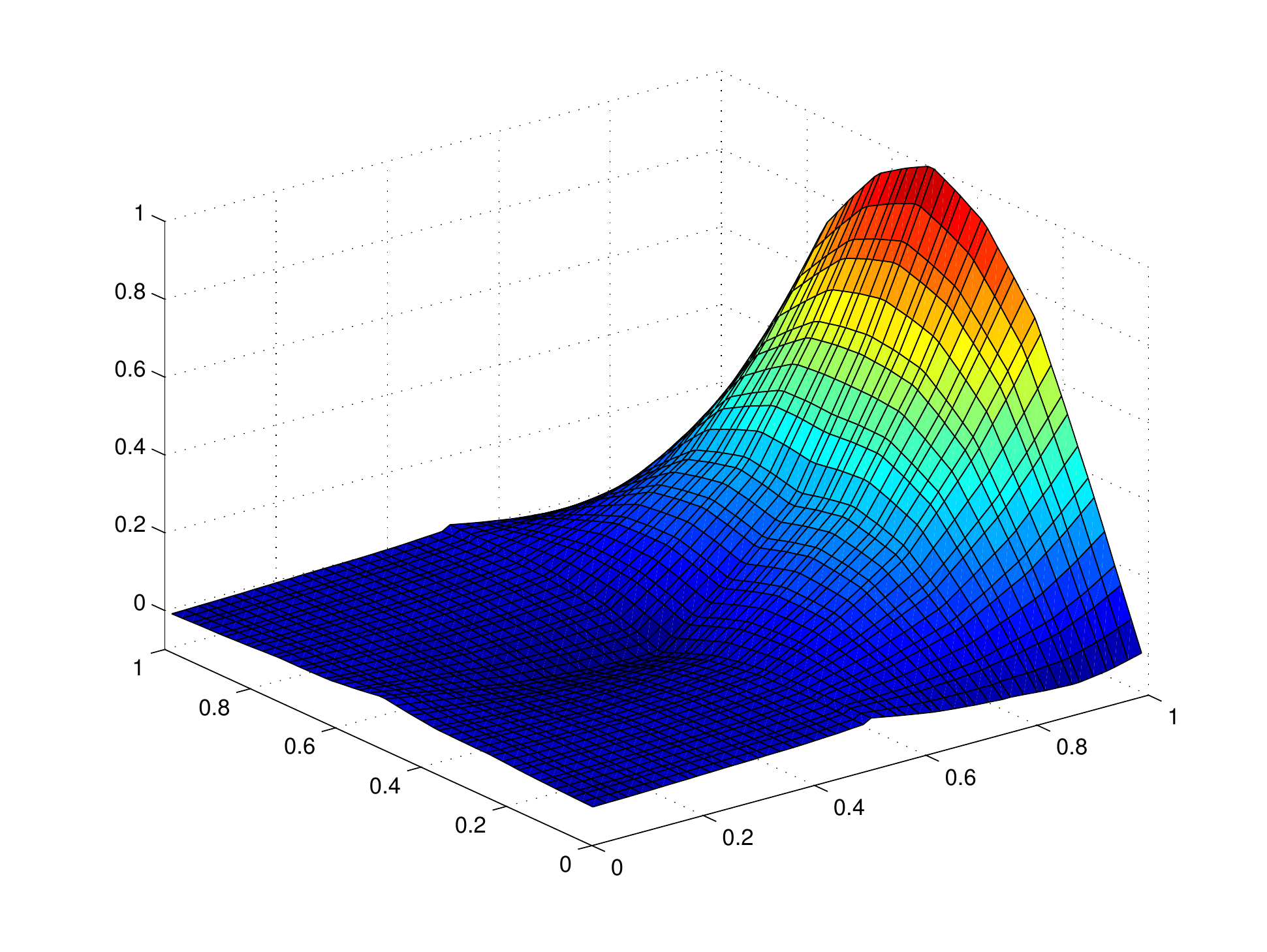}
(b) \includegraphics[width=.5\textwidth]{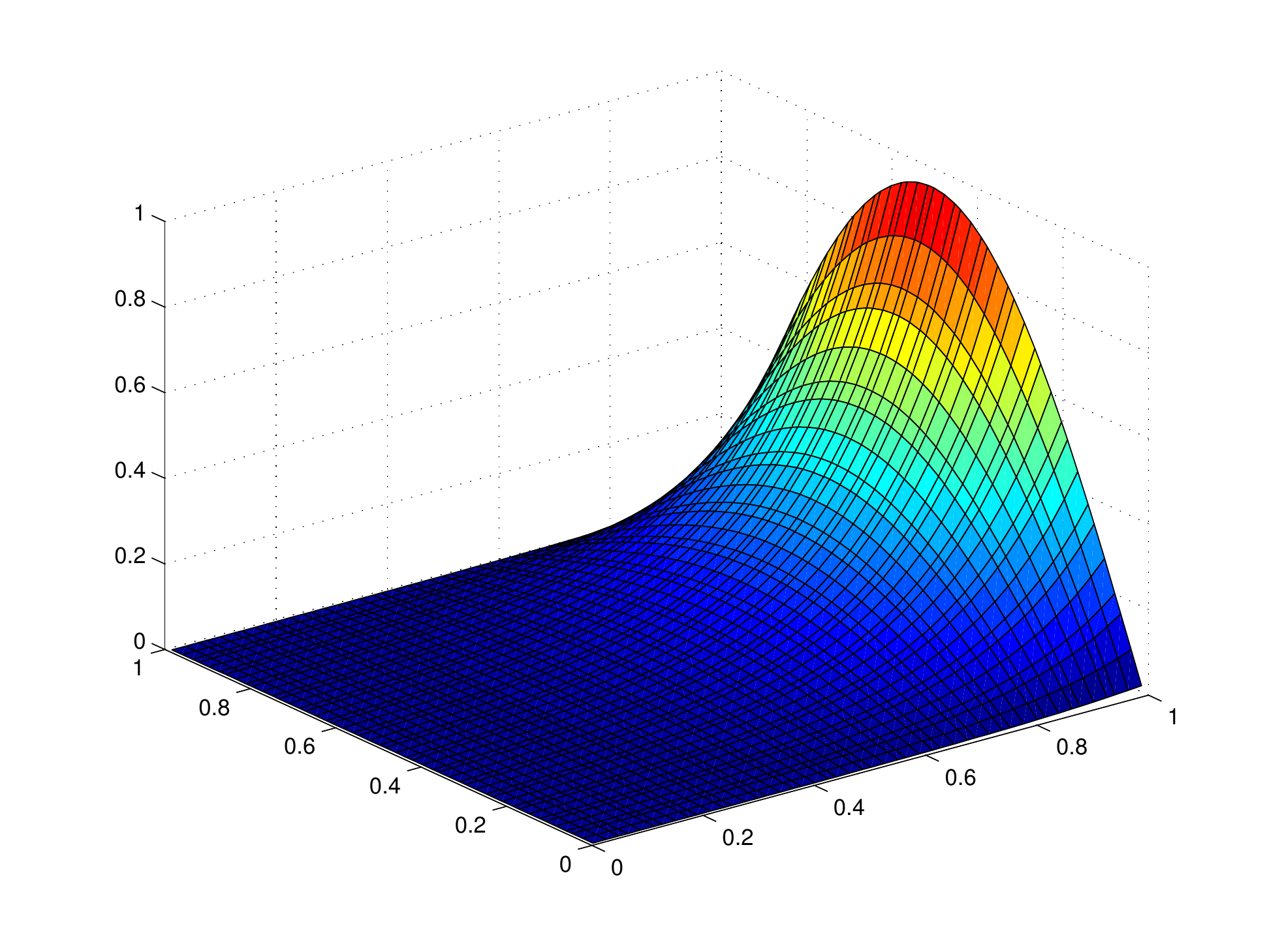}
\end{tabular}
  \caption{Example \ref{ex:const_4D}. (a) $k=1$ ; (b) $k=2$. $N=4$. Plotted along $x_1 = 0.4930, x_2 =0.4930$.
 }
 \label{fig:const_4D}
\end{figure}

\begin{table}
\caption{Numerical errors and orders of accuracy for Example \ref{ex:const_4D} computed by the space $\hat{\bV}_N^k$ with $ k=1, 2$ and indicated $N$.
}
\vspace{2 mm}
\centering
\begin{tabular}{c c c c c c c  c c }
\hline
N & $L^1$ error & order & $L^2$ error & order & $L^\infty$ error & order  &$H^1$ error & order \\
\hline

  &     & &  & $ k=1$    & &  &  & \\
\hline
 3   &2.44E-02  &          &  4.22E-02 &          &  3.31E-01 &            &3.91E-01 & \\
 4   &1.08E-02  & 1.18   & 2.08E-02 & 1.02   &  1.16E-01 &   1.51         &2.37E-01 &0.73 \\
 5   &3.68E-03  & 1.54   &  7.15E-03 &   1.54       &  9.33E-02 &  0.31          &1.22E-01 &0.96 \\
\hline
 &    & &  & $k=2$   &  &  & &   \\
 \hline
 2   &8.21E-04 &          &  1.34E-03 &          & 1.11E-02 &            &4.20E-02 & \\
 3   &1.76E-04  & 2.22  & 2.79E-04 & 2.27   &  2.76E-03 &   2.00         &1.20E-02 & 1.81 \\
 4   &3.32E-05  & 2.40   &  5.39E-05 & 2.37       & 8.76E-04 & 1.66        &3.18E-03 &1.91 \\
\hline

\end{tabular}
\label{table:4DconstBC}
\end{table}

\begin{table}
\caption{Sparsity and condition number of the stiffness matrix. Example \ref{ex:const_4D} computed by the space $\hat{\bV}_N^k$ with $ k=1, 2$. SGDOF is the number of degrees of freedom used for the sparse grid DG scheme. 
NNZ is the number of nonzero elements  in the stiffness matrix. Order=$\log(\textrm{NNZ})/\log(\textrm{SGDOF}).$}
\vspace{2 mm}
\centering
\begin{tabular}{c c c c c }
\hline
N &SGDOF & NNZ  & Order & Condition Number \\
\hline

&  & & $ k=1$   &\\
\hline
 3 &   1008  &12272 &   1.36  &4.27E+02 \\
4  &   3072  & 51712  &1.35  &2.26E+03   \\
5 & 8832 & 187008 & 1.34 & 9.27E+03  \\
\hline
 &  & & $k=2$     & \\
 \hline
2    & 1539  &19683    & 1.35 &7.40E+02    \\
3    & 5103 &102303 & 1.35 & 2.62E+03\\
4   &15552   &420336& 1.34 &9.72E+03 \\
\hline

\end{tabular}
\label{table:4Dconst_spy}
\end{table}

\begin{figure}[htp]
  \centering
  \begin{tabular}{c c}

(a)\includegraphics[width=.47\textwidth]{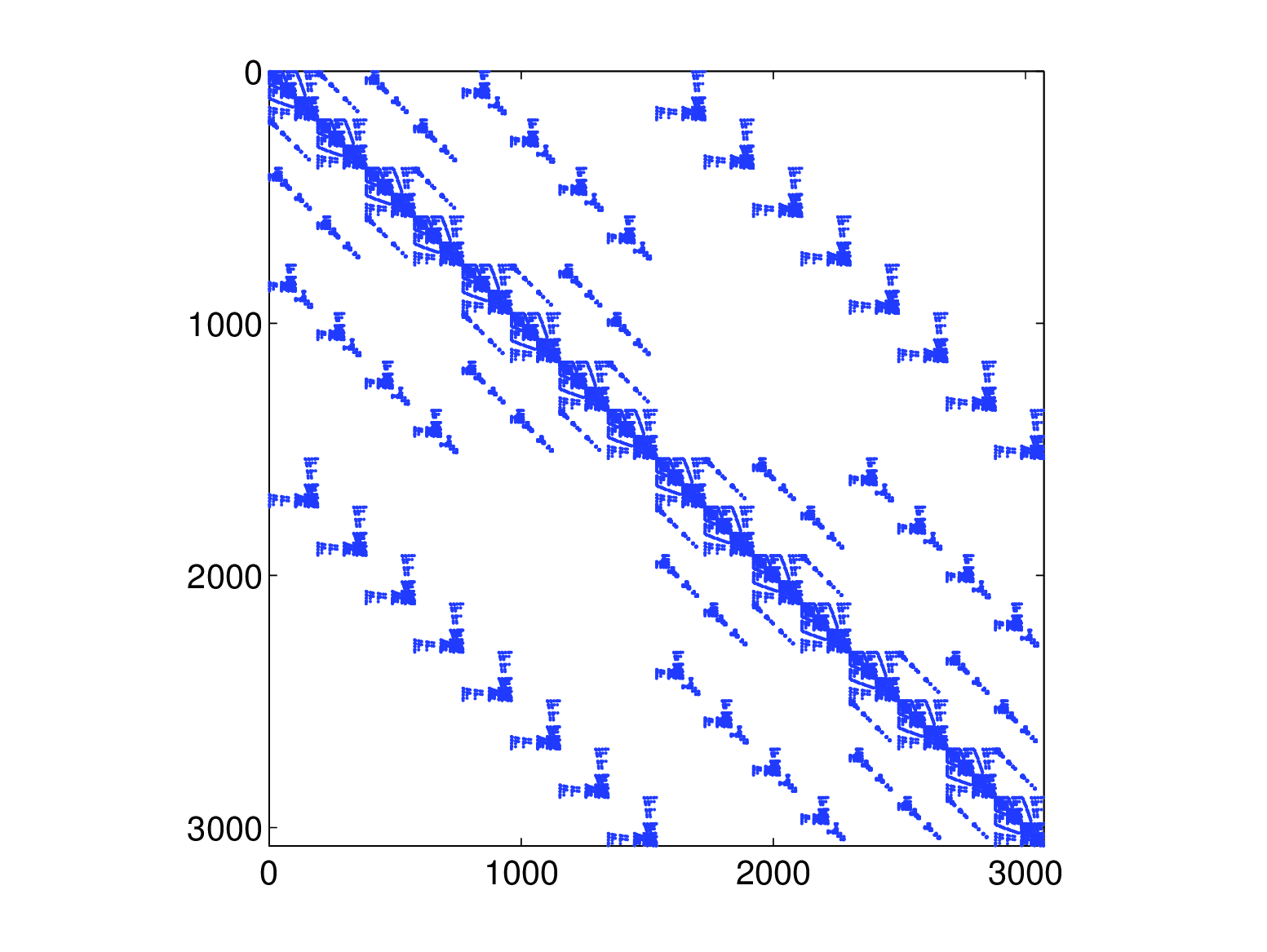}
(b)\includegraphics[width=.47\textwidth]{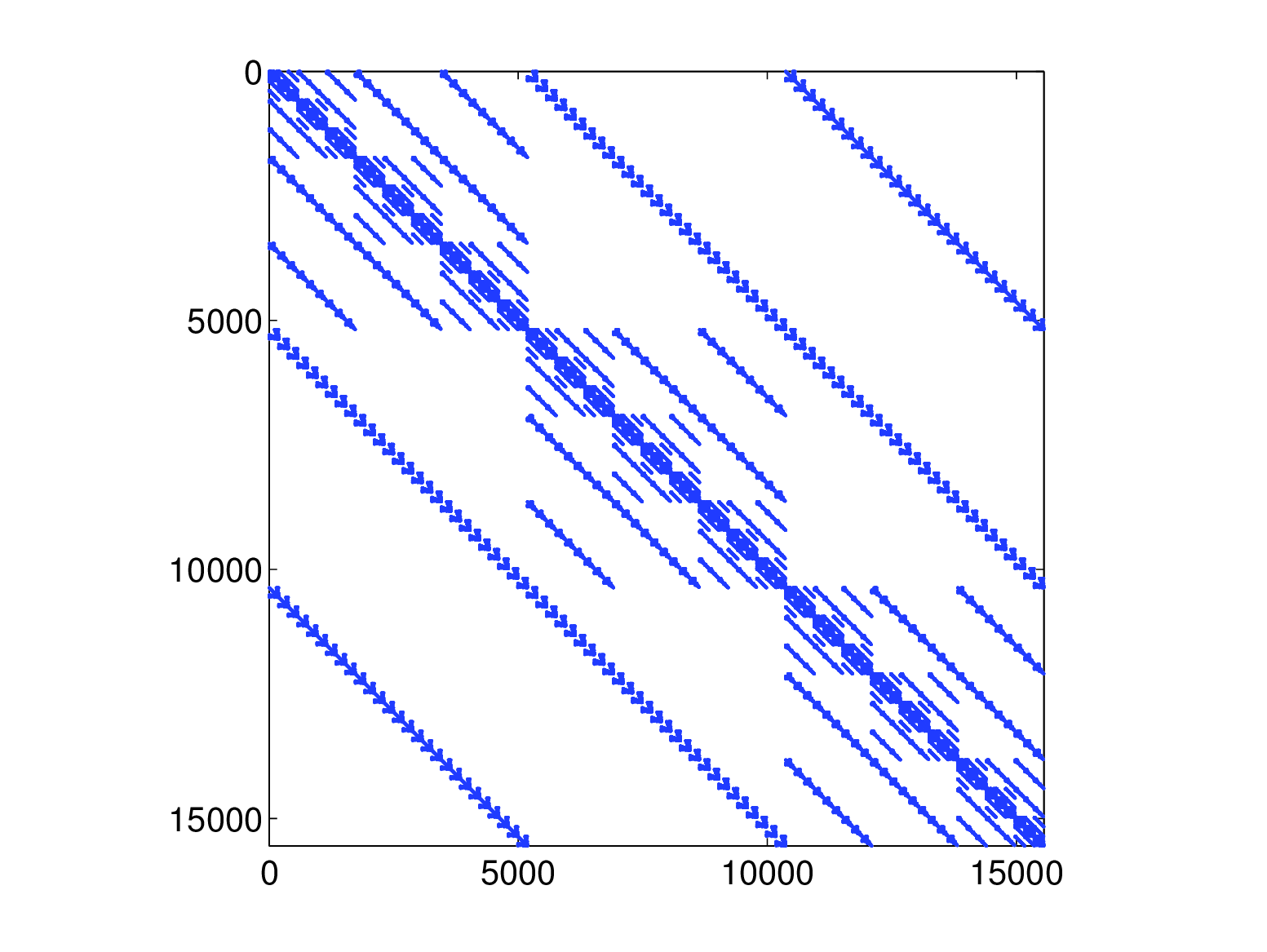}
\end{tabular}
  \caption{Example \ref{ex:const_4D}.  The sparsity patterns of the matrices computed by the space $\hat{\bV}_4^k$ with (a) $k=1$ and (b) $k=2$. Each dot represents a non-zero element in the stiffness matrix.}
 \label{fig:4Dconst_spy}
\end{figure}

\end{exa}

\newpage
\begin{exa}\rm We solve the following four-dimensional problem with smooth variable coefficient on $\Omega = [0,1]^4$.
\label{ex:4Dsmth}
\begin{align}
\begin{cases}
- \nabla \cdot ((\sin(x_1 x_2 x_3 x_4) +1 )\,  \nabla u) = f , &  \qquad \bx \in \Omega,\\
u = 0, &  \qquad \bx \in \partial\Omega,
  \end{cases}&
\end{align}
$f$ is a given function such that the exact solution is $u = \sin(\pi x_1)\, \sin(\pi x_2)\, \sin(\pi x_3)\, \sin(\pi x_4)$.
We provide the numerical results with $k=1$ and $k=2$  in Table \ref{table:4Dsmth}. The conclusion is similar to Example \ref{ex:const_4D}.

\begin{table} [!htp]
\caption{Numerical errors and   orders of accuracy for Example \ref{ex:4Dsmth}  computed by the space $\hat{\bV}_N^k$ with $ k=1, 2$ and indicated $N$. 
}
\centering
\begin{tabular}{c c c c c c c c c }
\hline
N &$L^1$ error & order & $L^2$ error & order & $L^\infty$ error & order  &$H^1$ error & order \\
\hline

  &   &  & &   $ k=1$    & &  &  & \\
\hline
3           &6.15E-02 &          &  8.97E-02 &          &  2.94E-01 &            &6.67E-01 & \\
4               &1.89E-02  & 1.70 &2.63E-02  &  1.77 & 2.54E-01  &0.21 & 3.20E-01& 1.06 \\
5            &  4.51E-03 &2.07 & 6.80E-03& 1.95 &  7.15E-02 &  1.83  &1.45E-01 & 1.14 \\
\hline
 &    & &   & $k=2$   & &  & &   \\
 \hline
2       &8.38E-04&         & 1.09E-03&          & 3.49E-03&          &3.74E-02&     \\
3       &1.62E-04  & 2.37 &   2.13E-04 & 2.36 &1.34E-03&1.38&  1.01E-02 & 1.90 \\
4    &2.97E-05  &2.44 & 3.91E-05 & 2.45 &3.80E-04 &  1.82 & 2.57E-03 & 1.97\\
\hline

\end{tabular}
\label{table:4Dsmth}
\end{table}

\end{exa}

\section{Conclusion}
\label{sec:Con}
In this paper, we develop a sparse grid IPDG method for elliptic equations on a box shaped domain. The method features a sparse finite element space which scales as $O(h^{-1}|\log_2 h|^{d-1})$ for $d$-dimensional problems,  translating into a significant cost reduction when $d$ is large. The traditional IPDG formulation can be readily used, and when combined with   approximation results for the new space, equip the scheme with  error estimate that is only slightly deteriorated in the energy norm. Numerical results  in up to four dimensions are shown to validate the performance of the method. Future work includes the study of adaptive algorithm for less regular solutions, extension of the method to more general domains, and developing sparse grid DG methods for other types of high-dimensional equations.

\bibliographystyle{abbrv}
\bibliography{ref_cheng}

\begin{thebibliography}{10}

\bibitem{achatz2003higher}
S.~Achatz.
\newblock Higher order sparse grid methods for elliptic partial differential
  equations with variable coefficients.
\newblock {\em Computing}, 71(1):1--15, 2003.

\bibitem{alpert1993class}
B.~Alpert.
\newblock A class of bases in {L}\^{}2 for the sparse representation of
  integral operators.
\newblock {\em SIAM J. Math. Anal.}, 24(1):246--262, 1993.

\bibitem{alpert2002adaptive}
B.~Alpert, G.~Beylkin, D.~Gines, and L.~Vozovoi.
\newblock Adaptive solution of partial differential equations in multiwavelet
  bases.
\newblock {\em J. Comput. Phys.}, 182(1):149--190, 2002.

\bibitem{archibald2011adaptive}
R.~Archibald, G.~Fann, and W.~Shelton.
\newblock Adaptive discontinuous {G}alerkin methods in multiwavelets bases.
\newblock {\em Appl. Numer. Math.}, 61(7):879--890, 2011.

\bibitem{arnold1982interior}
D.~Arnold.
\newblock An interior penalty finite element method with discontinuous
  elements.
\newblock {\em SIAM J. Numer. Anal.}, 19(4):742--760, 1982.

\bibitem{Arnold_2002_SIAM_DG}
D.~Arnold, F.~Brezzi, B.~Cockburn, and L.~Marini.
\newblock Unified analysis of discontinuous {Galerkin} methods for elliptic
  problems.
\newblock {\em SIAM J. Numer. Anal.}, 39:1749--1779, 2002.

\bibitem{babenko1960approximation}
K.~Babenko.
\newblock Approximation of periodic functions of many variables by
  trigonometric polynomials.
\newblock {\em Doklady Akademii Nauk SSSR}, 132(2):247--250, 1960.

\bibitem{baker1977finite}
G.~Baker.
\newblock Finite element methods for elliptic equations using nonconforming
  elements.
\newblock {\em Math. Comp.}, 31(137):45--59, 1977.

\bibitem{Bassi_1997_JCP_DFEM}
F.~Bassi and S.~Rebay.
\newblock A high-order accurate discontinuous finite element method for the
  numerical solution of the compressible {Navier-Stokes} equations.
\newblock {\em J. Comput. Phys.}, 131:267--279, 1997.

\bibitem{baszenski1992blending}
G.~Baszenski, F.-J. Delvos, and S.~Jester.
\newblock Blending approximations with sine functions.
\newblock In {\em Numerical Methods in Approximation Theory}, volume~9, pages
  1--19. Springer, 1992.

\bibitem{bellman1961adaptive}
R.~Bellman.
\newblock {\em Adaptive control processes: a guided tour}, volume~4.
\newblock Princeton University Press Princeton, 1961.

\bibitem{bungartz1997multigrid}
H.-J. Bungartz.
\newblock A multigrid algorithm for higher order finite elements on sparse
  grids.
\newblock {\em Electron T. Numer. Ana.}, 6:63--77, 1997.

\bibitem{bungartz1998finite}
H.-J. Bungartz.
\newblock {\em Finite elements of higher order on sparse grids}.
\newblock Shaker, 1998.

\bibitem{bungartz1998sparse}
H.-J. Bungartz and T.~Dornseifer.
\newblock Sparse grids: Recent developments for elliptic partial differential
  equations.
\newblock In {\em Multigrid Methods V}, pages 45--70. Springer, 1998.

\bibitem{bungartz2004sparse}
H.-J. Bungartz and M.~Griebel.
\newblock Sparse grids.
\newblock {\em Acta numer.}, 13:147--269, 2004.

\bibitem{calle2005wavelets}
J.~Calle, P.~Devloo, and S.~Gomes.
\newblock Wavelets and adaptive grids for the discontinuous {G}alerkin method.
\newblock {\em Numerical Algorithms}, 39(1-3):143--154, 2005.

\bibitem{castillo2000priori}
P.~Castillo, B.~Cockburn, I.~Perugia, and D.~Sch{\"o}tzau.
\newblock An a priori error analysis of the local discontinuous {G}alerkin
  method for elliptic problems.
\newblock {\em SIAM J. Numer. Anal.}, 38(5):1676--1706, 2000.

\bibitem{Ciarlet_1975_FEM_Elliptic}
P.~Ciarlet.
\newblock {\em The finite element method for elliptic problems}.
\newblock North-Holland, Amsterdam, 1975.

\bibitem{cockburn2001superconvergence}
B.~Cockburn, G.~Kanschat, I.~Perugia, and D.~Sch{\"o}tzau.
\newblock Superconvergence of the local discontinuous {G}alerkin method for
  elliptic problems on cartesian grids.
\newblock {\em SIAM J. Numer. Anal.}, 39(1):264--285, 2001.

\bibitem{Fengyan_Maxwell}
B.~Cockburn, F.~Li, and C.-W. Shu.
\newblock Locally divergence-free discontinuous {Galerkin} methods for the
  {Maxwell} equations.
\newblock {\em J. Comput. Phys.}, 194:588--610, 2004.

\bibitem{cockburn1998local}
B.~Cockburn and C.-W. Shu.
\newblock The local discontinuous {G}alerkin method for time-dependent
  convection-diffusion systems.
\newblock {\em SIAM J. Numer. Anal.}, 35(6):2440--2463, 1998.

\bibitem{dawson2004compatible}
C.~Dawson, S.~Sun, and M.~Wheeler.
\newblock Compatible algorithms for coupled flow and transport.
\newblock {\em Comput. Methods Appl. Mech. Eng.}, 193(23):2565--2580, 2004.

\bibitem{delvos1982d}
F.-J. Delvos.
\newblock $d$-variate {B}oolean interpolation.
\newblock {\em J. Approx. Theory.}, 34(2):99--114, 1982.

\bibitem{douglas1976interior}
J.~Douglas and T.~Dupont.
\newblock Interior penalty procedures for elliptic and parabolic galerkin
  methods.
\newblock In {\em Computing methods in applied sciences}, pages 207--216.
  Springer, 1976.

\bibitem{garcke2013sparse}
J.~Garcke and M.~Griebel.
\newblock {\em Sparse grids and applications}.
\newblock Springer, 2013.

\bibitem{gerhard2013adaptive}
N.~Gerhard and S.~M{\"u}ller.
\newblock Adaptive multiresolution discontinuous {G}alerkin schemes for
  conservation laws: multi-dimensional case.
\newblock {\em Comput. Appl. Math.}, pages 1--29, 2013.

\bibitem{gerstner1998numerical}
T.~Gerstner and M.~Griebel.
\newblock Numerical integration using sparse grids.
\newblock {\em Numerical algorithms}, 18(3-4):209--232, 1998.

\bibitem{gradinaru2007fourier}
V.~Gradinaru.
\newblock Fourier transform on sparse grids: code design and the time dependent
  {S}chr{\"o}dinger equation.
\newblock {\em Computing}, 80(1):1--22, 2007.

\bibitem{griebel1990parallelizable}
M.~Griebel.
\newblock A parallelizable and vectorizable multi-level algorithm on sparse
  grids.
\newblock In W.~Hackbusch, editor, {\em Parallel algorithms for partial
  differential equations}, volume~31 of {\em Notes on numerical fluid
  mechanics}, pages 94--100. 1991.

\bibitem{griebel1998adaptive}
M.~Griebel.
\newblock Adaptive sparse grid multilevel methods for elliptic {PDEs} based on
  finite differences.
\newblock {\em Computing}, 61(2):151--179, 1998.

\bibitem{griebel2005sparse}
M.~Griebel.
\newblock Sparse grids and related approximation schemes for higher dimensional
  problems.
\newblock In {\em Proceedings of the conference on foundations of computational
  mathematics}, Santander, Spain, 2005.

\bibitem{griebel2007sparse}
M.~Griebel and J.~Hamaekers.
\newblock Sparse grids for the {S}chr{\"o}dinger equation.
\newblock {\em ESAIM: Mathematical Modelling and Numerical Analysis},
  41(02):215--247, 2007.

\bibitem{griebel1995tensor}
M.~Griebel and P.~Oswald.
\newblock Tensor product type subspace splittings and multilevel iterative
  methods for anisotropic problems.
\newblock {\em Adv. Comput. Math.}, 4(1):171--206, 1995.

\bibitem{griebel1993multilevel}
M.~Griebel, C.~Zenger, and S.~Zimmer.
\newblock Multilevel {Gauss-Seidel-algorithms} for full and sparse grid
  problems.
\newblock {\em Computing}, 50(2):127--148, 1993.

\bibitem{griebel1999adaptive}
M.~Griebel and G.~Zumbusch.
\newblock Adaptive sparse grids for hyperbolic conservation laws.
\newblock In {\em Hyperbolic problems: theory, numerics, applications}, pages
  411--422. Springer, 1999.

\bibitem{haar}
A.~Haar.
\newblock Zur theorie der orthogonalen funktionensysteme.
\newblock {\em Mathematische Annalen}, 69(3):331--371, 1910.

\bibitem{hemker1995sparse}
P.~Hemker.
\newblock Sparse-grid finite-volume multigrid for 3{D}-problems.
\newblock {\em Adv. Comput. Math.}, 4(1):83--110, 1995.

\bibitem{hovhannisyan2014adaptive}
N.~Hovhannisyan, S.~M{\"u}ller, and R.~Sch{\"a}fer.
\newblock Adaptive multiresolution discontinuous {G}alerkin schemes for
  conservation laws.
\newblock {\em Math. Comp.}, 83(285):113--151, 2014.

\bibitem{iacono2011high}
F.~Iacono, G.~May, S.~M{\"u}ller, and R.~Sch{\"a}fer.
\newblock A high-order discontinuous {G}alerkin discretization with
  multiwavelet-based grid adaptation for compressible flows.
\newblock {\em AICES Preprint AICES-2011/08-1, RWTH Aachen}, 2011.

\bibitem{liem1995splitting}
C.~Liem, T.~L{\"u}, and T.~Shih.
\newblock {\em The splitting extrapolation method: a new technique in numerical
  solution of multidimensional problems}, volume~7.
\newblock World Scientific, 1995.

\bibitem{ma2009adaptive}
X.~Ma and N.~Zabaras.
\newblock An adaptive hierarchical sparse grid collocation algorithm for the
  solution of stochastic differential equations.
\newblock {\em J. Comput. Phys.}, 228(8):3084--3113, 2009.

\bibitem{mulder1989new}
W.~Mulder.
\newblock A new multigrid approach to convection problems.
\newblock {\em J. Comput. Phys.}, 83(2):303--323, 1989.

\bibitem{naik1993improved}
N.~Naik and J.~Van~Rosendale.
\newblock The improved robustness of multigrid elliptic solvers based on
  multiple semicoarsened grids.
\newblock {\em SIAM J. Numer. Anal.}, 30(1):215--229, 1993.

\bibitem{nobile2008sparse}
F.~Nobile, R.~Tempone, and C.~Webster.
\newblock A sparse grid stochastic collocation method for partial differential
  equations with random input data.
\newblock {\em SIAM J. Numer. Anal.}, 46(5):2309--2345, 2008.

\bibitem{oden1998discontinuoushpfinite}
J.~Oden, I.~Babu{\v{s}}ka, and C.~Baumann.
\newblock A discontinuous $hp$ finite element method for diffusion problems.
\newblock {\em J. Compt. Phys.}, 146(2):491--519, 1998.

\bibitem{pflaum1998multilevel}
C.~Pflaum.
\newblock A multilevel algorithm for the solution of second order elliptic
  differential equations on sparse grids.
\newblock {\em Numer. Math.}, 79(1):141--155, 1998.

\bibitem{riviere2001priori}
B.~Rivi{\`e}re, M.~Wheeler, and V.~Girault.
\newblock A priori error estimates for finite element methods based on
  discontinuous approximation spaces for elliptic problems.
\newblock {\em SIAM J. Numer. Anal.}, 39(3):902--931, 2001.

\bibitem{schwab2008sparse}
C.~Schwab, E.~S{\"u}li, and R.~Todor.
\newblock Sparse finite element approximation of high-dimensional
  transport-dominated diffusion problems.
\newblock {\em ESAIM: Mathematical Modelling and Numerical Analysis},
  42(05):777--819, 2008.

\bibitem{shen2010sparse}
J.~Shen and L.-L. Wang.
\newblock Sparse spectral approximations of high-dimensional problems based on
  hyperbolic cross.
\newblock {\em SIAM J. Numer. Anal.}, 48(3):1087--1109, 2010.

\bibitem{shen2010efficient}
J.~Shen and H.~Yu.
\newblock Efficient spectral sparse grid methods and applications to
  high-dimensional elliptic problems.
\newblock {\em SIAM J. Sci. Comput.}, 32(6):3228--3250, 2010.

\bibitem{smolyak1963quadrature}
S.~Smolyak.
\newblock Quadrature and interpolation formulas for tensor products of certain
  classes of functions.
\newblock In {\em Dokl. Akad. Nauk SSSR}, volume~4, pages 240--243, 1963.

\bibitem{temlyakov1986approximations}
V.~Temlyakov.
\newblock Approximations of functions with bounded mixed derivative.
\newblock {\em Trudy Matematicheskogo Instituta im. VA Steklova}, 178:3--113,
  1986.

\bibitem{vuik2014multiwavelet}
M.~Vuik and J.~Ryan.
\newblock Multiwavelet troubled-cell indicator for discontinuity detection of
  discontinuous {G}alerkin schemes.
\newblock {\em J. Comput. Phys.}, 270:138--160, 2014.

\bibitem{wang2009wkb}
W.~Wang and C.-W. Shu.
\newblock The {WKB} local discontinuous {G}alerkin method for the simulation of
  {S}chr{\"o}dinger equation in a resonant tunneling diode.
\newblock {\em J. Sci. Comput.}, 40(1-3):360--374, 2009.

\bibitem{zixuanthesis}
Z.~Wang.
\newblock {\em Discontinuous {Galerkin} methods for {Hamilton-Jacobi} equations
  and high-dimensional elliptic equations}.
\newblock PhD thesis, Michigan State University, 2015.

\bibitem{wheeler1978elliptic}
M.~Wheeler.
\newblock An elliptic collocation-finite element method with interior
  penalties.
\newblock {\em SIAM J. Numer. Anal.}, 15(1):152--161, 1978.

\bibitem{xiu2007efficient}
D.~Xiu.
\newblock Efficient collocational approach for parametric uncertainty analysis.
\newblock {\em Comm. Comput. Phys.}, 2(2):293--309, 2007.

\bibitem{xiu2005high}
D.~Xiu and J.~Hesthaven.
\newblock High-order collocation methods for differential equations with random
  inputs.
\newblock {\em SIAM J. Sci. Comput.}, 27(3):1118--1139, 2005.

\bibitem{yuan2006discontinuous}
L.~Yuan and C.-W. Shu.
\newblock Discontinuous {G}alerkin method based on non-polynomial approximation
  spaces.
\newblock {\em J. Comput. Phys.}, 218(1):295--323, 2006.

\bibitem{zenger1991sparse}
C.~Zenger.
\newblock Sparse grids.
\newblock In {\em Parallel Algorithms for Partial Differential Equations,
  Proceedings of the Sixth GAMM-Seminar}, volume~31, 1990.

\end{thebibliography}

\end{document}